\documentclass{paper}
\usepackage{graphicx}
\usepackage{graphicx,amssymb,amsmath,amsfonts,mathtext,color,wrapfig}

\def\be{\begin{eqnarray}}
\def\ee{\end{eqnarray}}

\newcommand{\eq}[1]{\begin{equation} #1 \end{equation}}

\hoffset= - 0.5in         

\begin{document}

\thispagestyle{empty}

\baselineskip14pt

\hfill ITEP/TH-10/12

\bigskip

\bigskip

\centerline{\LARGE{On Undulation Invariants of Plane Curves
}}

\vspace{5ex}

\centerline{\large{\emph{A.Popolitov\footnote{Institute for Theoretical and Experimental Physics, Moscow, Russia; popolit@itep.ru} and Sh.Shakirov\footnote{Department of Mathematics, UC Berkeley, USA and Center for Theoretical Physics, UC Berkeley, USA; \newline Institute for Theoretical and Experimental Physics, Moscow, Russia; shakirov@math.berkeley.edu}}}}

\vspace{4ex}

\centerline{ABSTRACT}

\bigskip

{\footnotesize
One of the general problems in algebraic geometry is to determine algorithmically whether or not a given geometric object, defined by explicit polynomial equations (e.g. a curve or a surface), satisfies a given property (e.g. has singularities or other distinctive features of interest). A classical example of such a problem, described by A.Cayley and G.Salmon in 1852, is to determine whether or not a given plane curve of degree $r > 3$ has undulation points -- the points where the tangent line meets the curve with multiplicity four. They proved that there exists an invariant of degree $6(r - 3)(3 r - 2)$ that vanishes if and only if the curve has undulation points. In this paper we give explicit formulae for this invariant in the case of quartics ($r=4$) and quintics ($r=5$), expressing it as the determinant of a matrix with polynomial entries, of sizes $21 \times 21$ and $36 \times 36$ respectively.
}

\section{Introduction}

It is widely accepted that \emph{linear algebra}, that deals with objects such as systems of linear equations, matrices, determinants, eigenvalues and eigenvectors, is of exceptional importance not only as a branch of pure mathematics, but mostly as a practical tool and language for almost any calculation with several indeterminates. Because linear algebra is developed so well, we have a nearly complete understanding of various linear phenomena in nature, ranging from the most fundamental (e.g. superpositions of wave functions in quantum mechanics) to the most practical (e.g. distributions of electricity currents in circuits).

This situation stands in great constrast with our understanding of non-linear phenomena (from confinement in atomic nuclei to turbulence in water flows) which is far less complete. Such phenomena are usually governed by non-linear differential equations, more complicated than the linear ones. One possible solution, often used in practice, is to try to consider a non-linear process or object as approximately linear. This leads at best to approximate solutions, what may be satisfactory in some cases but gives very limited understanding and in general is not enough. For this reason, in our days development of exact methods to deal with non-linear problems attracts more and more attention.

From a mathematician's standpoint, the most fundamental way to attack the problem of understanding non-linearity would be to think first about non-linear algebra, i.e. about  generalization of linear algebra to the non-linear setting, with linear equations replaced by generic algebraic equations, their zero sets -- lines and planes -- replaced by algebraic curves and surfaces, matrices (two-index arrays) replaced by tensors (multi-index arrays), and so on. Such generalization,  started by the prominent scientists of the 19th century, is today known as \emph{algebraic geometry}.

\pagebreak

On the physical side, the same endeavor towards better understanding of non-linear problems has influenced many fields, including theory of dynamical systems, integrability, topological field theory, and eventually string theory. These and other related branches of theoretical physics, unified by their motivation to preserve and adequately describe the non-linear mathematical structure of physical problems, are often described as \emph{mathematical physics}.

It may seem that the non-linear algebraic structures are too simple to capture the sophisticated essense of the non-linear phenomena in nature. However, this is not quite true. Despite realistic problems are described by differential, not algebraic equations, it is well-known that many of the puzzling non-linear phenomena -- such as chaotic dynamics or fractals -- appear already in relatively simple algebraic problems \cite{Mandelbrot}. Vice versa, physical methods and ideas often lead to further progress in algebraic geometry: a fine example is solution of Riemann-Schottky problem with the help of KP equation and integrability theory \cite{RiemannSchottky}. For these reasons, many believe that collaboration between algebraic geometry and mathematical physics is beneficial for both sides. A research program that approaches non-linear algebra from this point of view has been recently suggested in \cite{Nolinal}.

Non-linear algebraic equations are interesting and have many applications, but finding their exact solutions is typically hard. If a non-linear equation is encountered in practice, it is most often treated with one or another type of approximation. Exact solutions are rarely considered. One of the achievments of algebraic geometry is that it allows to extract exact information about the solutions without actually computing the solutions theirselves. This is done by taking look at the quantities called \emph{algebraic invariants} of the system.

The simplest example of an invariant is just the determinant in linear algebra -- a condition of solvability for a system of $n$ homogeneous linear equations in $n$ unknowns. Such a system is solvable (has a non-zero solution $x_1, \ldots, x_n$) if and only if a certain polynomial in coefficients, called the determinant, vanishes: $\det = 0$. For example, for a pair of linear equations
$$
\left\{
\begin{array}{l}
f_{1} x_1 + f_{2} x_2 = 0 \\
\\
g_{1} x_1 + g_{2} x_2 = 0
\end{array}
\right.
$$
the determinant equals
$$
\det = f_1 g_2 - f_2 g_1
$$
and the system is solvable if and only if $\det = 0$. The determinant is invariant of the action of the group $SL(n)$, that changes the basis in the underlying linear space while preserving the volume. This is clear since the fact of existence of solutions is independent on the choice of basis.

Algebraic geometry provides a direct non-linear generalization of determinant, called \emph{resultant} -- a condition of solvability for a system of $n$ homogeneous algebraic equations in $n$ unknowns. Such a system is solvable (has a non-zero solution $x_1, \ldots, x_n$) if and only if a certain polynomial in coefficients, called the resultant, vanishes: $R = 0$. For example, for a pair of quadratic equations
$$
\left\{
\begin{array}{l}
f_{11} x_1^2 + f_{12} x_1 x_2 + f_{22} x_2^2 = 0 \\
\\
g_{11} x_1^2 + g_{12} x_1 x_2 + g_{22} x_2^2 = 0
\end{array}
\right.
$$
the resultant equals
$$
R = f_{11}^2 g_{22}^2 - 2 g_{22} f_{11} f_{22} g_{11}-g_{12} g_{22} f_{12} f_{11} + f_{11} f_{22} g_{12}^2 -g_{12} f_{12} f_{22} g_{11}+g_{11} g_{22} f_{12}^2 +f_{22}^2 g_{11}^2
$$
and the system is solvable if and only if $R = 0$. In other words, knowledge of the resultant gives a direct way to find out, whether the non-linear system has a solution or not -- something that otherwise is not straightforward to achieve.

By the same argument as the determinant, the resultant is an invariant of $SL(n)$. The difference is, however, that in linear algebra there is essentially only one invariant -- the determinant -- while in the non-linear case there are many. All of these invariants carry different pieces of information about the solutions of the system. The resultant, for example, tells us whether the system has at least one solution at all (zero solution not counted) while other invariants can distinguish the cases when the system has more than one solution, or certain prescribed number of solutions, or solutions with certain prescribed characteristics, etc.

This shows that invariants allow to learn a lot about the solutions without actually finding the solutions -- a very useful property. Already for this reason explicit construction of invariants, as well as study of their general properties, is an important direction of work. Many problems in this field were stated and solved by the great mathematicians of the 19th century: Cayley, Bezout, Klein, Sylvester and others made contributions to this classical area of research. Still, in many cases their results have the form of existence theorems, not explicit formulae.

The reason for this is that invariants are usually very large polynomials, so computations with them, except for the simplest cases, have relatively high computational complexity, certainly too high for calculations "by hand". In a sense, algebraic geometry of the 19th century lacked the computational power of modern computers. It is today that construction of invariants and study of their properties goes to a new level, and many problems left open in the 19th century can now be solved. One such problem is a topic of this paper.

\section{The undulation problem}

\begin{wrapfigure}{r}{0.4\textwidth}
\vspace{-4ex}
  \begin{center}
    \includegraphics[width=0.4\textwidth]{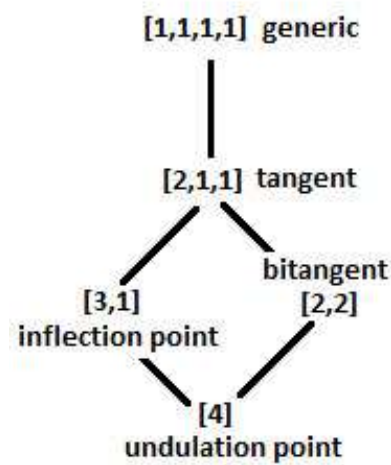}
  \end{center}
  \vspace{-1ex}
  \caption{Types of intersection of a plane curve and a line. The case $r = 4$.}
  \vspace{1ex}
\end{wrapfigure}
This paper is devoted to one particular problem where invariants turn out to be useful, due to A.Cayley and G.Salmon \cite{CayleySalmon}. To explain the problem we first start with informal discussion and then proceed to the rigorous formulation. Let us consider the projective plane ${\mathbb P}^2$ with homogeneous coordinates $x_1:x_2:x_3$, and a plane curve
$$
P(x_1,x_2,x_3) = 0
$$
defined by a homogeneous irreducible polynomial
$$
P(x_1,x_2,x_3) = \sum\limits_{i+j+k=r} C_{ijk} \ x_1^i \ x_2^j \ x_3^k
$$
of degree $r$. For $r = 1$ the curve is called a line, for $r = 2$ a quadric, for $r = 3$ a cubic, for $r = 4$ a quartic, and so on. For brevity, let $C$ denote the set of coefficients $C_{ijk}$. This set consists of $\# C = (r+1)(r+2)/2$ coefficients.

By Bezout theorem, any line in ${\mathbb P}^2$ crosses this curve in exactly $r$ points, if counted with multiplicities. For this reason, the types of possible intersections can be put into correspondence with partitions of the degree, $r = m_1 + m_2 + \ldots$, where the parts $m_i$ of the partition are the multiplicities of intersection points. Clearly, if a line is generic it intersects the curve in $r$ distinct points with all multiplicities 1, i.e. it corresponds to the partition $(1,1,\ldots,1)$. The simplest non-generic intersection occurs for the \emph{tangent line} to a curve: then one of the intersection points has multiplicity 2 (the point of tangency) while all the other intersection points have multiplicity 1. This type of intersection corresponds to partition $(2,1,1, \ldots)$.

\begin{figure}[t]
  \begin{center}
    \includegraphics[width=0.2\textwidth]{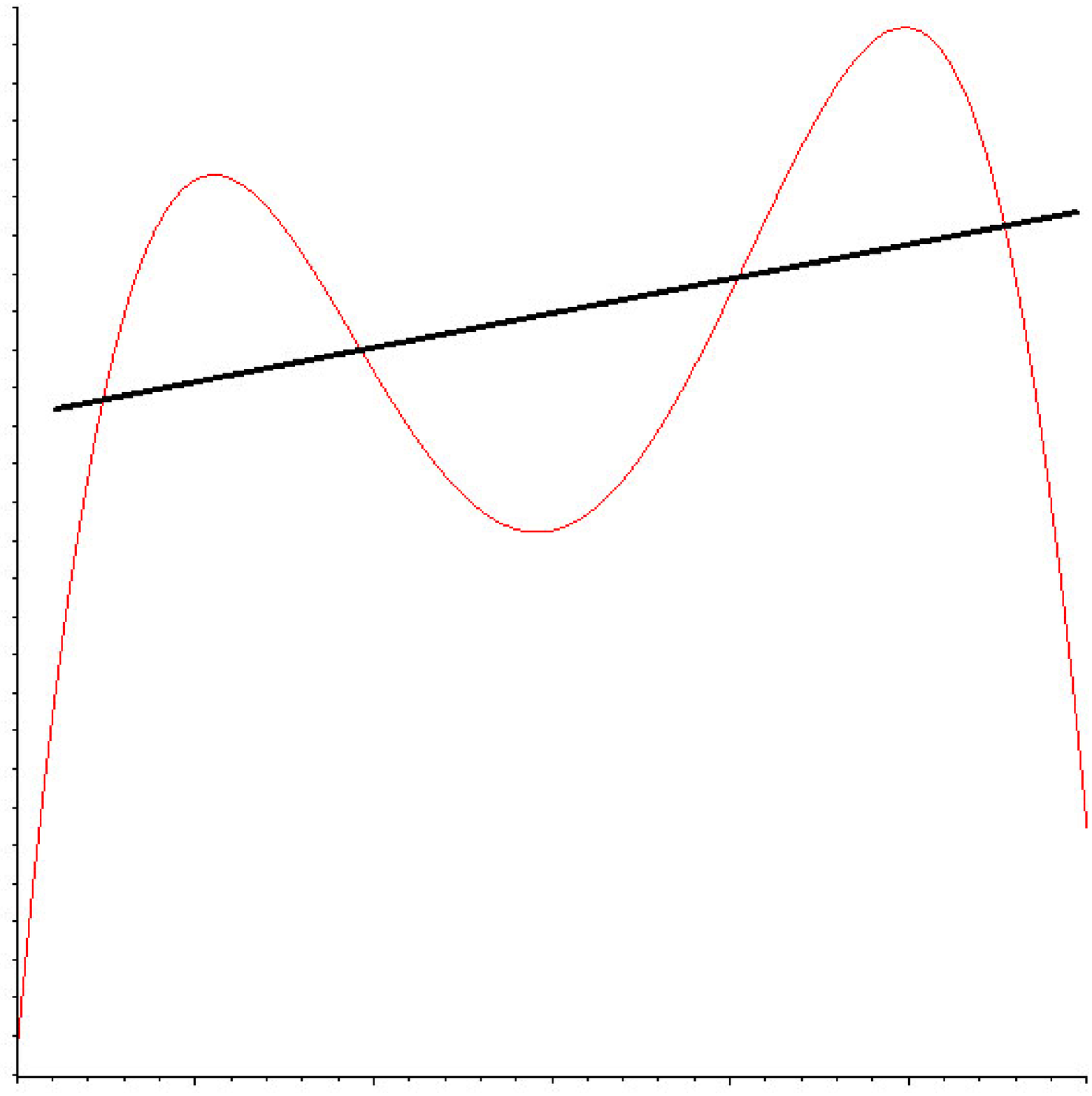}
\ \ \ \ \ \ \
    \includegraphics[width=0.2\textwidth]{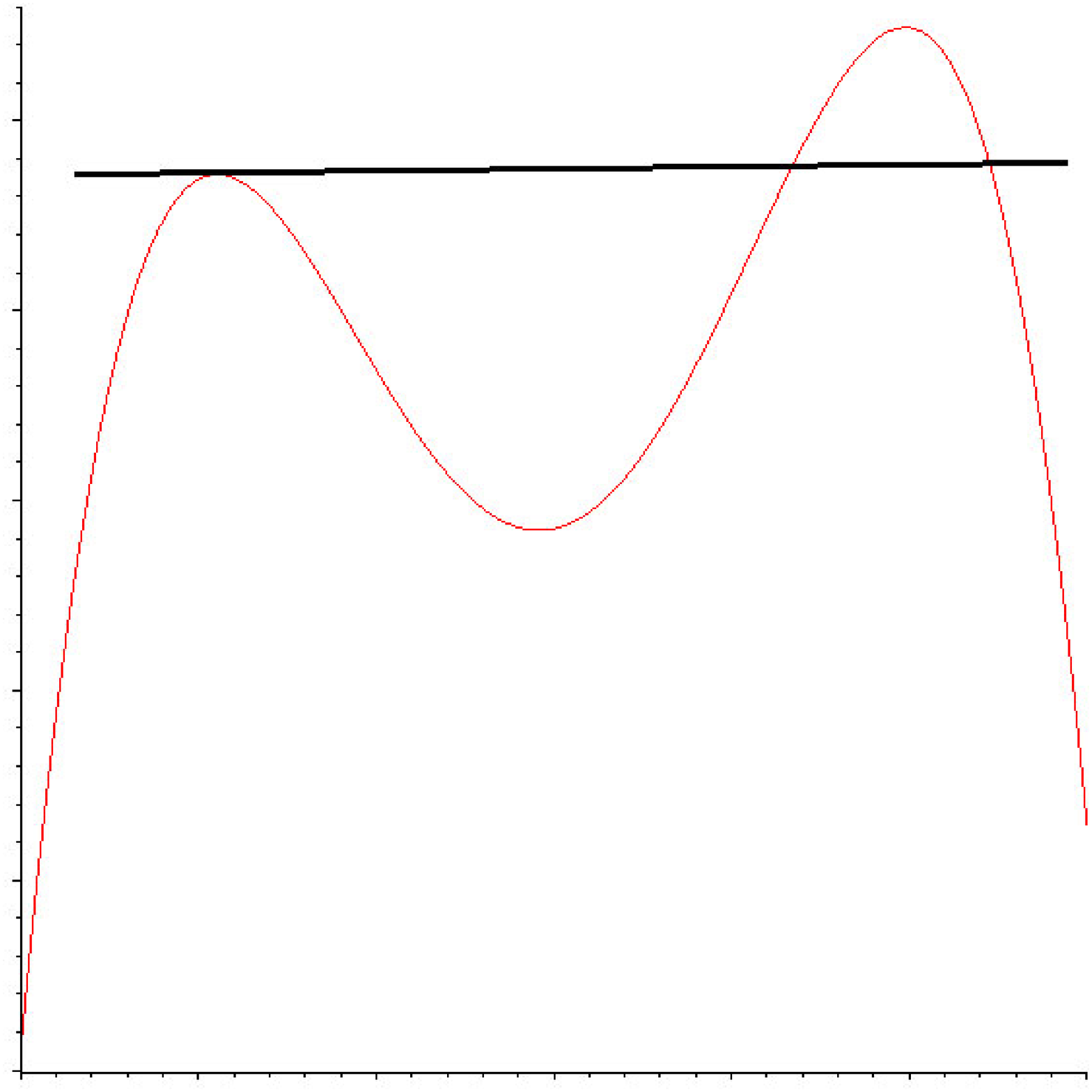}
\ \ \ \ \ \ \
    \includegraphics[width=0.2\textwidth]{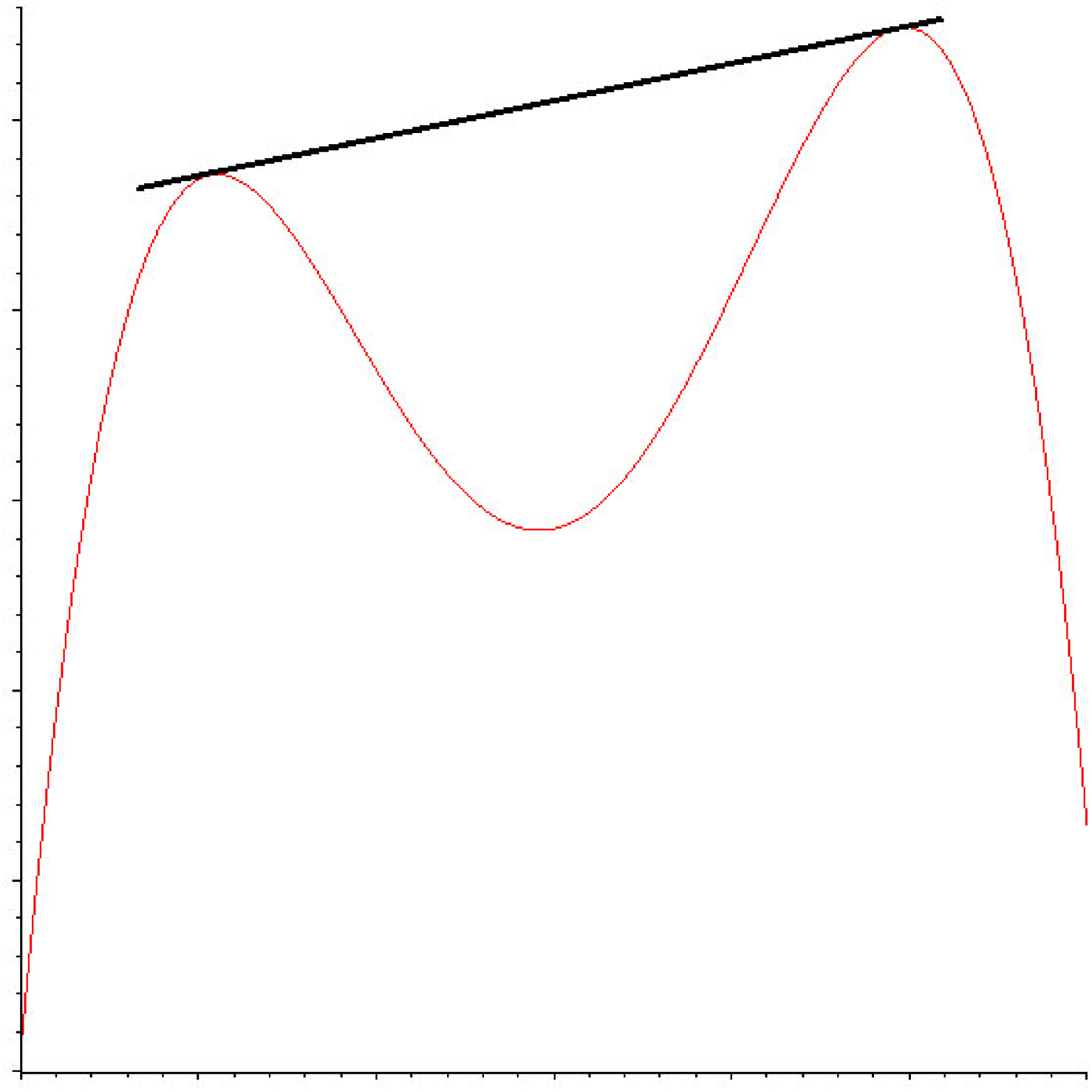}
\ \ \ \ \ \ \
    \includegraphics[width=0.2\textwidth]{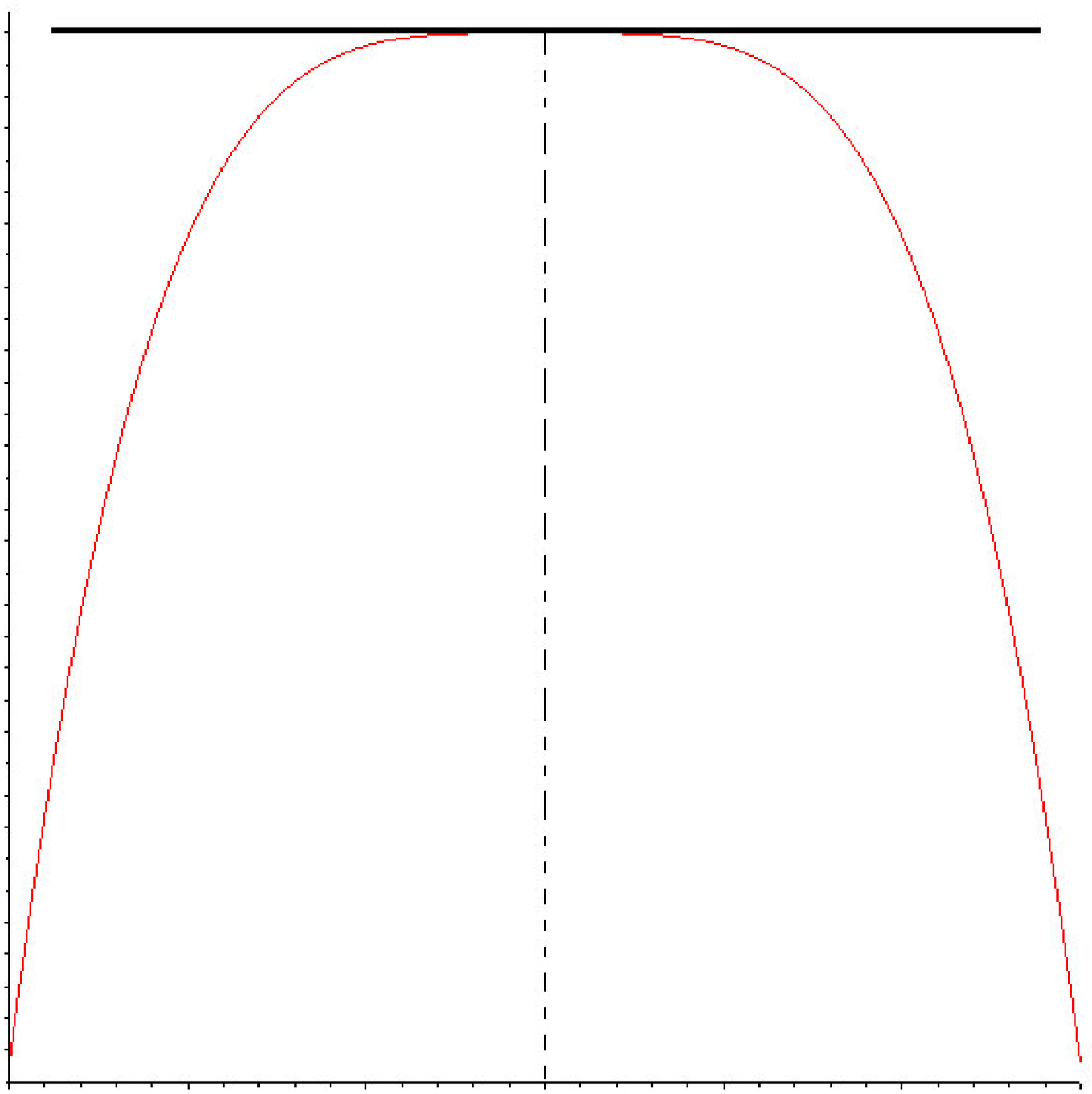}
  \end{center}
  \caption{Four intersections of a line and a plane quartic of types $(1,1,1,1), (2,1,1), (2,2)$ and $(4)$, respectively. Pictures in real section. The undulation point is marked by a dashed line. }
\end{figure}

The next-to-simplest types of intersection are, respectively, $(3,1,\ldots, 1)$ and $(2,2,1,\ldots,1)$; the former situation is called a \emph{line of inflection}, while the latter - a \emph{bitangent}, since in this case the line is sumultaneously tangent to a curve in two distinct points. One can continue further by considering lines of type $(4,1,\ldots,1)$, $(3,2,1,\ldots,1)$, and so on. These generally do not have given names, with one notable exception: a line of type $(4,1,\ldots,1)$ is called a \emph{line of undulation}, and the corresponding point of intersection is called an \emph{undulation point}. Note that conventionally intersection with at least one part of multiplicity higher than 4 is also called an undulation point:

\paragraph{2.1. Definition.}A point is called an undulation point of a plane curve $P(x_1,x_2,x_3) = 0$ if a tangent line at that point meets the curve with multiplicity four or higher.

\paragraph{} It is classically known (and it is easy to estimate from degree counting) that a given plane curve has only finitely many bitangent lines. Moreover, in generic position all of these bitangent lines are distinct. This immediately implies that a generic plane curve has no undulation points -- since undulation points can be thought of as a result of merging of two bitangents. For a plane curve to posess undulation points, it should be non-generic, i.e. $P$ should satisfy some algebraic equation(s). Cayley and Salmon studied this question and proved the following existence theorem:

\paragraph{2.2. Theorem \cite{CayleySalmon}}There exists a unique up to rescaling function $I$, which is a homogeneous polynomial in $(r+1)(r+2)/2$ coefficients $C$ of degree $6(r-3)(3r-2)$, such that
$$
I(C) = 0
$$
is a necessary and sufficient condition for the curve $P(x_1,x_2,x_3) = 0$ to have undulation points.

\paragraph{}This important statement illustrates once again the usefulness of invariant theory in practical considerations. For a given plane curve, a direct method of investigating whether this curve has undulation points or not is almost hopeless: it requires solving a system of non-linear algebraic equations that is, for generic $r$, not even solvable in radicals. If, on the other hand, one knows the polynomial $I$ explicitly, it is immediate to decide, are there any undulation points or not. It goes without saying that the polynomial $I$ is invariant under the action of $SL(3)$, by the usual reasoning (existence of undulation points is a condition that holds irrespective of the choice of basis).

However, this theorem only justifies existence of such an invariant. In practice, it is useful to have not only that but also an explicit formula, that can be used to explicitly compute the invariant. By the \textbf{undulation problem} we mean the problem of finding an explicit polynomial formula for the undulation invariant $I$, that could be realized as an algorithm and used in practice. Such a formula was not given in \cite{CayleySalmon}, and the aim of this paper is to fill this gap. We will show that $I(C)$ is given by a determinant with polynomial entries.

We believe that one of the reasons that the explicit formula that we present here was not found before is the complexity of the invariants: even in the simplest non-trivial case of quartic curves, for $r = 4$, the invariant $I$ is a homogeneous polynomial of degree $60$ in the $15$ variables $ C = \{ C_{400}, C_{310}, C_{301}, C_{220}, C_{211}, C_{202}, C_{130}, $ $ C_{121}, C_{112}, C_{103}, C_{040}, C_{031}, C_{022}, C_{013}, C_{004} \} $ and has

$$
\dfrac{(15+60-1)!}{(15-1)!60!} =
456002537343216 \simeq 4.5 \cdot 10^{14}
$$
\smallskip\\
monomials! This terrifying "growth of length" is a very typical phenomenon in invariant theory and, generally speaking, in non-linear algebra (see, e.g. the appendix to the book version of \cite{Nolinal} for another example). It is, however, not terrifying at all if appropriate structures and key properties of the object are identified: surprisingly or not, these enormous invariant polynomials often posess a lot of nice properties and can be expressed by simple, elegant formulae. In the above example, the invariant of degree $60$ turns out to be a determinant of a relatively simple $21 \times 21$ matrix. To find and explore such formulae, it is often beneficial to use modern computers and software (e.g. MAPLE, see s. 7 for more detailed discussion of technical tools we used).

\section{Algebraic reformulation}

The undulation problem is so far defined in the geometric terms of tangent lines and multiplicities. To proceed to solution of this problem, the following algebraic reformulation will be useful:

\paragraph{3.1. Lemma.}A line

\begin{align}
v(x_1,x_2,x_3) = v_1 x_1 + v_2 x_2 + v_3 x_3 = 0
\end{align}
\smallskip\\
is an undulation line of a plane curve $P(x_1,x_2,x_3) = 0$, iff $P$ can be decomposed in a form

\begin{align}
P(x_1,x_2,x_3) = u(x_1,x_2,x_3)^4 h(x_1,x_2,x_3) + v(x_1,x_2,x_3) w(x_1,x_2,x_3)
\label{Decomp}
\end{align}
\smallskip\\
where $u(x_1,x_2,x_3) = u_1 x_1 + u_2 x_2 + u_3 x_3$ is some linear polynomial, and $h(x_1,x_2,x_3), w(x_1,x_2,x_3)$ are some homogeneous polynomials of degrees $r-4$ and $r-1$ respectively.

\paragraph{Proof $(\Rightarrow)$} Suppose that $P$ can be written in the form (\ref{Decomp}), and let $(X_1 : X_2 : X_3)$ be the (unique, because of irreducibility of $P$) intersection point of the two linear polynomials

$$
\left\{
\begin{array}{ll}
u_1 X_1 + u_2 X_2 + u_3 X_3 = 0 \\
\\
v_1 X_1 + v_2 X_2 + v_3 X_3 = 0
\end{array}
\right.
$$
namely,

$$
(X_1 : X_2 : X_3) = \big( u_2 v_3 - u_3 v_2 : \ u_3 v_1 - u_1 v_3 : \ u_1 v_2 - u_2 v_3 \big)
$$
\smallskip\\
Straightforward calculation then gives

$$
P \Big|_{x=X} = 0, \ \ \ \ \dfrac{\partial P}{\partial x_i} \Big|_{x=X} = v_i
$$
\smallskip\\
i.e. the line $v(x_1,x_2,x_3) = 0$ is a tangent line to the curve $P(x_1,x_2,x_3) = 0$ at the point $(X_1 : X_2 : X_3)$. Calculating the intersection multiplicity using the axioms that define it, we find

$$
I_X( P, v ) = I_X( u^4 h + v w, v) = I_X ( u^4 h, v) = 4 I_X(u, v) + I_X(h, v) \geq 4
$$
\smallskip\\
i.e. $X$ is indeed an undulation point, \textbf{q.e.d.}

\paragraph{Proof $(\Leftarrow)$}  Suppose that $v(x_1,x_2,x_3) = 0$ is a tangent line to $P(x_1,x_2,x_3) = 0$ at $(X_1 : X_2 : X_3)$ with intersection multiplicity at least 4. Without loss of generality (using an appropriate change of coordinates on ${\mathbb C}^3$) we can assume that $(X_1 : X_2 : X_3) = (0:0:1)$ and $v(x_1,x_2,x_3) = x_1$. It is well-known that under these assumptions the intersection multiplicity of a line $v(x_1,x_2,x_3) = 0$ and a curve $P(x_1,x_2,x_3) = 0$ is nothing but the multiplicity of zero as a root of a single-variable polynomial $p(z) = P(0,z,1)$. In our case the multiplicity is at least 4, i.e. $p(z)$ has a form

$$
P(0,z,1) = z^4 w(z)
$$
\smallskip\\
for some polynomial $w(z)$. Passing to homogeneous polynomials, we find

$$
P(0,x_2,x_3) = x_2^4 w(x_2,x_3)
$$
for some homogeneous polynomial $w(x_2,x_3)$, thus

$$
P(x_1,x_2,x_3) = x_2^4 w(x_2,x_3) + x_1 h(x_1,x_2,x_3)
$$
\smallskip\\
for some homogeneous polynomial $h(x_1,x_2,x_3)$. Therefore, $P$ indeed has the form (\ref{Decomp}), \textbf{q.e.d.}

\pagebreak

\section{The undulation ideal}

To solve the undulation problem, the following algebraic object will be of key importance.

\paragraph{4.1. Definition.} The \emph{undulation ideal} ${\cal I}$ is the set of all polynomials in the variables $C,v_1,v_2,v_3$ that vanish whenever $v_1 x_1 + v_2 x_2 + v_3 x_3 = 0$ is the undulation line of the curve $P(x_1,x_2,x_3) = 0$:
$$
{\cal I} = \big\{ \ f \in {\mathbb C}[C,v_1,v_2,v_3] \ \big| \ f(C,v_1,v_2,v_3) = 0 \ \mbox{ if } v_1 x_1 + v_2 x_2 + v_3 x_3 = 0 \big.
$$
$$
\big. \mbox{ is an undulation line for the curve } P(x_1,x_2,x_3) = 0 \ \big\} $$

\paragraph{} The motivation to consider such an ideal essentially comes from the Cayley-Salmon theorem 2.2. In algebraic language, this theorem can be stated as a fact that the simpler ideal
$$
{\cal I}^{\prime} = \big\{ \ f \in {\mathbb C}[C] \ \big| \ f(C) = 0 \ \mbox{ if } P(x_1,x_2,x_3) = 0 \mbox{ has at least one undulation line } \big\}
$$
is generated by a single element -- the undulation invariant:
$$
{\cal I}^{\prime} = \big< I(C) \big>
$$
Following the general wisdom \emph{"to understand something, deform/generalize it"} we propose to extend ${\cal I}^{\prime}$ to a bigger ideal ${\cal I}$, in a hope that this could reveal an additional structure and thus shed some light on the nature of the element $I(C)$. As we will see, this will happen to be the case.

The ideal ${\cal I}$ admits three useful gradings. The first two are the obvious gradings w.r.t. the total degree in all the coefficients $C$, and the total degree in all the coefficients $v$:

$$
\deg_C: {\cal I} \rightarrow {\mathbb Z}_+, \ \ \ \deg_C(f) = \sum\limits_{i+j+k=r} \deg_{C_{ijk}}(f)
$$
$$
\deg_v: {\cal I} \rightarrow {\mathbb Z}_+, \ \ \ \deg_v(f) = \deg_{v_1}(f) + \deg_{v_2}(f) + \deg_{v_3}(f)
$$
\smallskip\\
The last, third, grading is more refined. In a certain sense, it keeps track of degrees in the coordinate variables $x_1,x_2,x_3$. Speaking more precisely, let us define
$$
{\overline \deg}: {\cal I} \rightarrow {\mathbb Z}^3_+
$$
to be the unique homomorphism of rings such that

$$
\overline{\deg}(v_1) = (1,0,0), \ \ \ \overline{\deg}(v_2) = (0,1,0), \ \ \ \overline{\deg}(v_3) = (0,0,1), \ \ \ \overline{\deg}(C_{ijk}) = (i,j,k)
$$
\smallskip\\
With respect to these gradings, the ideal ${\cal I}$ decomposes into a direct sum of graded components.

\paragraph{4.2. Definition.} Let ${\cal I}_{n,m}$ be the graded components of ${\cal I}$ w.r.t. the first two gradings, and ${\cal I}_{n,m_1,m_2,m_3}$ be the graded components w.r.t. all the three gradings:

$$
{\cal I}_{n,m} = \big\{ f \in {\cal I} \big| \deg_C(f) = n, \ \deg_v(f) = m \big\}
$$

$$
{\cal I}_{n,m_1,m_2,m_3} = \big\{ f \in {\cal I} \big| \deg_C(f) = n, \ \overline{\deg}(f) = (m_1,m_2,m_3) \big\}
$$
\smallskip\\
Note, that the $v$-degree of ${\cal I}_{n,m_1,m_2,m_3}$ is $m_1 + m_2 + m_3 - n$. Note also, that the Cayley-Salmon ideal is nothing but the graded component of ${\cal I}$ that has zero $v$-degree, ${\cal I}^{\prime} = \oplus_n {\cal I}_{n,0}$.

To understand the structure of the graded components ${\cal I}_{n,m}$, various methods can be used. It is interesting that for the purpose of this paper it is enough to use the most basic and straightforward approach possible: direct calculation of the spaces ${\cal I}_{n,m}$ using only Lemma 3.1. and linear algebra. This direct approach can be summarized in two statements.

\paragraph{4.3. Corollary.} Each ${\cal I}_{n,m}$ is a finite-dimensional linear space.

\paragraph{4.4. Corollary.} Each ${\cal I}_{n,m}$ can be computed as a solution to a finite linear system of equations.

\paragraph{Proof.}Corollary 4.3. follows directly from the definition. For Corollary 4.4., denote

$$
s_{ijk}(u,h,v,w) = \mbox{ coefficient in front of } x_1^i x_2^j x_3^k \mbox{ in } \big( u^4 h + v w \big)
$$
\smallskip\\
with $u,h,v,w$ as in Lemma 3.1. Then, a homogeneous polynomial $f \in {\mathbb C}[C,v_1,v_2,v_3]$ of degrees $\deg_C(f) = n, \deg_v(f) = m$ belongs to ${\cal I}_{n,m}$ iff the following system of equations is satisfied:

\eq{
\label{poly_systems}
f\big( s(u,h,v,w), v_1, v_2, v_3 \big) = 0 \ \ \ \forall u,v,h,w
}
\smallskip\\
This is a system of finitely many linear equations, where the coefficients of the polynomial $f$ are treated as indeterminates. For given pair of natural numbers $n,m$, there are only finitely many these coefficients. Therefore, for any given $n,m$, one can (at least in principle) write and explicitly solve the corresponding linear system, obtaining ${\cal I}_{n,m}$ as its solution space. \textbf{Q.e.d.}

\paragraph{} Despite the size and complexity of the above linear systems grows quite fast with $n,m$, we will see below that this straightforward approach suffices to investigate the simplest properties of the ideal ${\cal I}$. In particular, in the next section we will use this approach to find several lowest ${\cal I}_{n,m}$ for $r = 4$ and show, that the elements of these linear spaces can be naturally put together to form a $21 \times 21$ matrix, determinant of which is the Cayley-Salmon invariant. This is the main new result of current paper, that calls for further research in the nearby directions.Then, in the next following section, we will do the same for $r = 5$ and obtain similar results, thus giving evidence that the $r = 4$ result is not an accident but rather the first step towards generalizations.

\pagebreak

\emph{}

\section{Determinantal formula for $r = 4$}

Using the approach explained in the previous section, we obtain the following

\paragraph{Theorem 5.1.} For plane quartics ($r = 4$) the dimensions dim ${\cal I}_{n,m}$ of a few lowest graded components of the undulation ideal are given by the following numbers:

\[
\begin{array}{c|ccccccccccccccccccccc}
n \backslash m & 0 & 1 & 2 & 3 & 4 & 5 & 6 & 7 & \ldots \\
\hline
0 & 0 & 0 & 0 & 0 & 0 & 0 & 0 & 0 & \ldots \\
1 & 0 & 0 & 0 & 0 & 0 & 0 & 0 & 0 & \ldots \\
2 & 0 & 0 & 0 & 0 & 1 & \textbf{3} & 21 & 45 & \ldots \\
3 & 0 & 0 & 0 & 0 & 15 & \textbf{63} & 325 & \ldots & \ldots \\
\ldots & \ldots & \ldots & \ldots & \ldots & \ldots & \ldots & \ldots & \ldots & \ldots
\end{array}
\]
\paragraph{Proof.} Direct calculation via Corollaries 4.3., 4.4.

\paragraph{} The spaces ${\cal I}_{2,5}$ and ${\cal I}_{3,5}$ are spanned, as linear spaces, by 3 and 63 polynomials, resp.:
$$
{\cal I}_{2,5} = \big\{ \sum\limits_{i=1}^{3} c_i \alpha_i \big| c_1,\ldots,c_3 \in {\mathbb C} \big\}
$$
$$
{\cal I}_{3,5} = \big\{ \sum\limits_{i=1}^{63} c_{i} \beta_{i} \big| c_1,\ldots,c_{63} \in {\mathbb C} \big\}
$$
Since ${\cal I}$ is an ideal, a product of any element of ${\cal I}_{2,5}$ and any element of $C$ belongs to ${\cal I}_{3,5}$. This implies that ${\cal I}_{3,5}$ is a direct sum of a 45-dimensional subspace spanned by such products, and the complementary 18-dimensional subspace. Let $\beta_1, \ldots, \beta_{18}$ be the basis elements of that 18-dimensional subspace. Together with the 3 basis elements of ${\cal I}_{2,5}$, they form a set of 21 linearly independent polynomials of degree 5 in $v_1,v_2,v_3$. At the same time, the dimension of the space of homogeneous polynomials of degree 5 in three variables $v_1,v_2,v_3$ is exactly 21! This allows us to arrange these 3 + 18 polynomials into a $21 \times 21$ matrix, with the following remarkable property.

\paragraph{Theorem 5.2.} Let ${\cal M}$ be the $21 \times 21$ matrix, the rows of which are obtained by expanding the 21 polynomials $\alpha_1, \ldots, \alpha_3; \beta_1, \ldots, \beta_{18}$ in the 21 homogeneous monomials of degree 5 in $v_1,v_2,v_3$. Then the determinant of this matrix is the Cayley-Salmon undulation invariant of plane quartics:
\begin{equation}
\addtolength{\fboxsep}{5pt}
\boxed{
\begin{gathered}
I(C)_{r=4} = \det\limits_{21 \times 21} {\cal M}
\end{gathered}
}\label{DetermQuartic}
\end{equation}

\paragraph{Proof.} By construction, if the curve $P(x_1,x_2,x_3) = 0$ posesses an undulation line $V_1x_1 + V_2 x_2 + V_3x_3 = 0$, then all the polynomials $\alpha_1, \ldots, \alpha_3; \beta_1, \ldots, \beta_{18}$ vanish at $v = V$. This implies that the 21-dimensional vector, components of which are the 21 monomials of degree 5 in $V_1,V_2,V_3$, belongs to the kernel of ${\cal M}$. Hence ${\cal M}$ is degenerate \pagebreak whenever the curve $P(x_1,x_2,x_3) = 0$ posesses an undulation line. Hence, its determinant is an element of the Cayley-Salmon undulation ideal:
$$
\det\limits_{21 \times 21} {\cal M} \in {\cal I}^{\prime} = \big< I(C) \big>
$$
By Cayley-Salmon theorem, this ideal is generated by a unique element -- the Cayley-Salmon undulation invariant $I(C)$ -- and therefore det ${\cal M}$ has to be proportional to $I(C)$:
$$
\det\limits_{21 \times 21} {\cal M} = i(C) \cdot I(C)
$$
Here $i(C)$ is some polynomial in $C$. It is easy to calculate its degree:
$$
\deg_C i(C) = \deg_C \det\limits_{21 \times 21} {\cal M} - \deg_C I(C) = 3 \cdot 2 + 18 \cdot 3 - 60 = 0
$$
So $i(C)$ does not depend on $C$, i.e. is just a constant. The determinant of ${\cal M}$ thus has the same degree as the Cayley-Salmon undulation invariant, and coincides with it up to an overall constant that can be always put to 1 (since $I(C)$ is itself defined up to rescaling). \textbf{Q.e.d.}

\section{Representation theory structure}

\begin{wrapfigure}{r}{0.25\textwidth}
\vspace{-10ex}
  \begin{center}
    \includegraphics[width=0.25\textwidth]{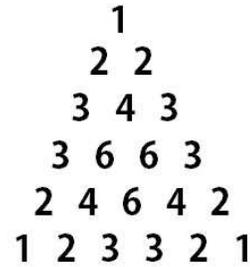}
  \end{center}
  \vspace{-4ex}
  \caption{Structure of ${\cal I}_{3,5}$.}
  \vspace{1ex}
\end{wrapfigure}
As usual, the symmetry group of the problem (in our case, $SL(3)$) acts on the space of solutions, decomposing it into irreducible representations. To find this decomposition, the easiest way is to consider, instead of the graded components ${\cal I}_{n,m}$, the more refined components ${\cal I}_{n,m_1,m_2,m_3}$. In complete analogy, their dimensions and spanning polynomials can be computed via Corollaries 4.3., 4.4. For $n = 3, m = 5$ this gives the following triangle of integers, as shown on Fig. 3. Decomposing this triangle into the usual multiplicity diagrams of irreducible representations of $SL(3)$, we find
\begin{center}
\includegraphics[width=0.5\textwidth]{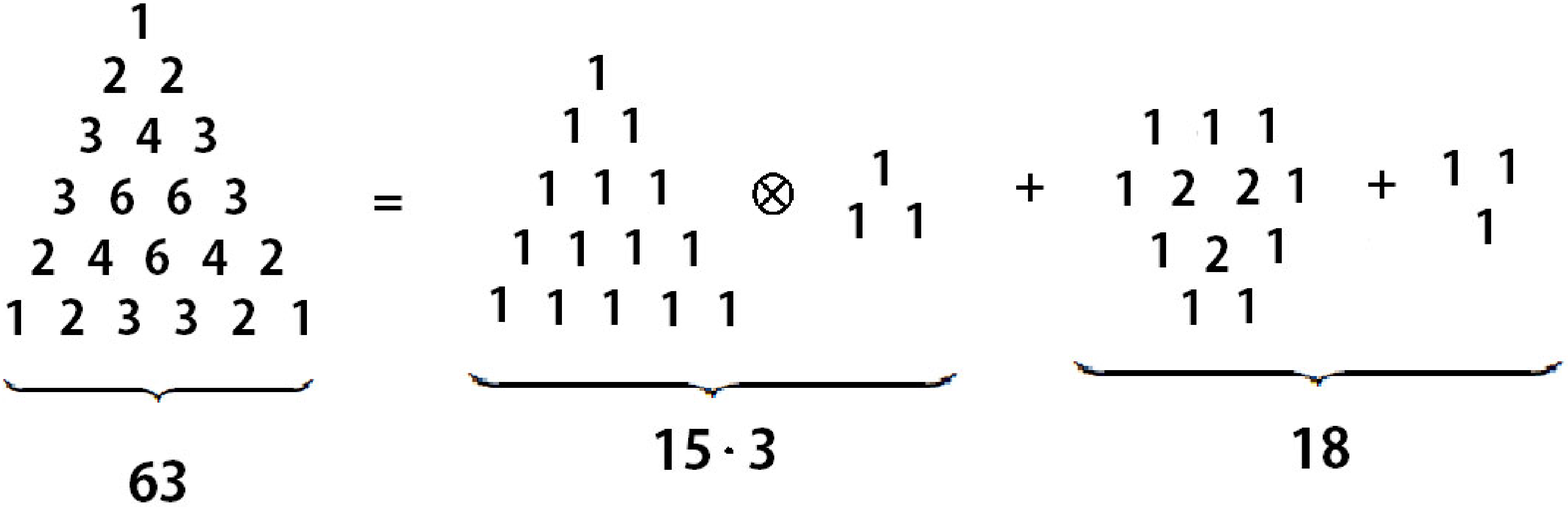}
\end{center}
This indicates that the 63 basis polynomials of ${\cal I}_{3,5}$ consist of: the non-interesting $15 \cdot 3 = 45$ polynomials obtained from ${\cal I}_{3,4}$, the 15 polynomials transforming in the irreducible representation (3,2), and the 3 polynomials transforming in the irreducible representation $(1,1) = ({\bar 1})$. This completes the description of the decomposition of ${\cal I}_{3,5}$ into irreducibles. Finally, it is an easy exercise to check that ${\cal I}_{2,5}$ is itself an irreducible representation of $SL(3)$ -- namely, just the fundamental representation (1). Together, these $21 = 15+3+3$ polynomials, transforming in representations $(3,2) \oplus (1,1) \oplus (1)$, form the matrix ${\cal M}$. Explicit expressions for these polynomials, necessary for the actual computations, are given in the Appendix.

\section{Determinantal formula for $r = 5$}

Having solved the undulation problem for $r = 4$, it is natural to go further and consider the case $r = 5$ -- that of quintics. For plane quintics, the undulation invariant has degree $6(r-3)(3r-2) = 156$. Despite this is an impressively large degree, in fact almost nothing changes compared to the case of plane quartics, and we are able to present the following main Theorems 7.1 and 7.2.

\paragraph{Theorem 7.1.} For plane quintics ($r = 5$) the dimensions dim ${\cal I}_{n,m}$ of a few lowest graded components of the undulation ideal are given by the following numbers:

\[
\begin{array}{c|ccccccccccccccccccccc}
n \backslash m & \leq 5 & 6 & 7 & \ldots \\
\hline
0 & 0 & 0 & 0 & \ldots \\
1 & 0 & 0 & 0 & \ldots \\
2 & 0 & 6 & \textbf{15} & \ldots \\
3 & 0 & 126 & 315 & \ldots \\
\\
\ldots & \ldots & \ldots & \ldots & \ldots \\
\\
6 & 0 & 63756 & \textbf{159411} & \ldots \\
\\
\end{array}
\]
\paragraph{Proof.} Direct calculation via Corollaries 4.3., 4.4. Note that, despite theoretically this calculation is just as direct as for Theorem 5.1, in practice it is significantly harder because sizes of linear systems \ref{poly_systems}, which define generators of the ideal, become so large, that solving them with MAPLE (and even finding their rank) is no longer possible. We tackle this technical problem by utilizing the ``linbox'' linear algebra package. However, we are not using it directly, writing program to compute the ranks in pure C. Rather we use SAGE, which is a great tool for mathematicians, that is written in Python and has bindings for ``linbox''. To glue our MAPLE and SAGE code together (that is to convert linear systems from polynomial form of notation, tgenerated by the former to sparse-matrix form, understood by the latter) we use simple Perl script. We believe, that such a pattern of using several distinct computational and modeling tools, each of which is well suited for a particular task -- rather than using all-in-one swiss-knives -- and then glueing them together with help of scripting languages (such as Perl, Python and Lisp) will become more and more common in mathematical physics.

However, fair compuration of ranks (and then the dimensions of solution spaces, which we are interested in) over field of zero characteristic is still impossible in this case, even with help of all these techniques. So we compute them over field integers modulo 6361 (which, we believe, is a sufficiently large prime number).

Having said this, we are able to obtain the following table of dimensions of solution spaces for ${\cal I}_{6,m_1,m_2,m_3}$ (horizontal axis is $m_1$ and vertical - $m_2$, while $m_3 = 7 - m_2 - m_1$).

\begin{center}
\includegraphics[width=0.8\textwidth]{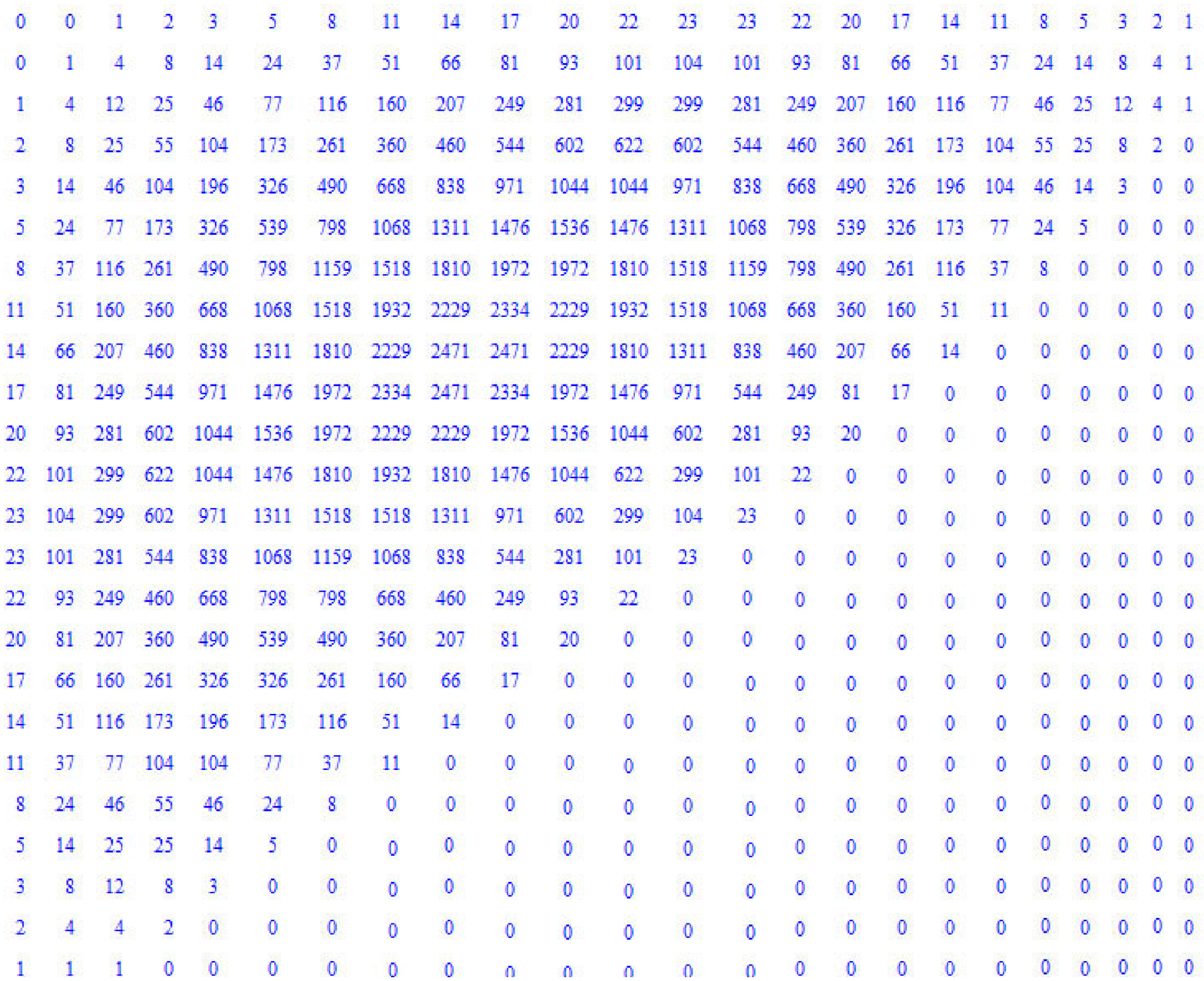}
\end{center}

Thus, the spaces ${\cal I}_{2,7}$ and ${\cal I}_{6,7}$ are spanned, as linear spaces, by 15 and 159411 polynomials, resp. As before, let us denote the basis polynomials in these spaces by $\alpha_1, \ldots, \alpha_{15}$ and $\beta_1, \ldots, \beta_{159411}$. Since ${\cal I}$ is an ideal, a product of any element of ${\cal I}_{2,7}$ and any polynomial of degree 4 of $C$ belongs to ${\cal I}_{6,7}$. This implies that ${\cal I}_{6,7}$ is a direct sum of a subspace spanned by such products, and the complementary subspace of dimension
$$
159411 - 15 \cdot {\rm dim} \Big( S^4 \ {\rm span}(C) \Big) = 159411 - 15 \cdot \dfrac{(21 + 4 - 1)!}{4!(21-1)!} = 159411 - 159390 = 21
$$
Let $\beta_1, \ldots, \beta_{21}$ be the basis elements of that 21-dimensional subspace. Together with the 15 basis elements of ${\cal I}_{2,6}$, they form a set of 36 linearly independent polynomials of degree 7 in $v_1,v_2,v_3$. At the same time, the dimension of the space of homogeneous polynomials of degree 7 in three variables $v_1,v_2,v_3$ is exactly 36! This allows us to arrange these 15 + 21 polynomials into a $36 \times 36$ matrix, with the following remarkable property.

\paragraph{Theorem 7.2.} Let ${\cal M}$ be the $36 \times 36$ matrix, the rows of which are obtained by expanding the 36 polynomials $\alpha_1, \ldots, \alpha_{15}; \beta_1, \ldots, \beta_{21}$ in the 36 homogeneous monomials of degree 7 in $v_1,v_2,v_3$. Then the determinant of this matrix is the Cayley-Salmon undulation invariant of plane quintics:
\begin{equation}
\addtolength{\fboxsep}{5pt}
\boxed{
\begin{gathered}
I(C)_{r=5} = \det\limits_{36 \times 36} {\cal M}
\end{gathered}
}\label{DetermQuintic}
\end{equation}

\paragraph{Proof.} Analogously to Theorem 5.2, simply follows from the fact that $\det {\cal M}$ vanishes when the curve has undulation points, and has the correct degree: $15 \cdot 2 + 21 \cdot 3 = 156$.

\section{Conclusion}

In this paper we solved the undulation problem for plane quartics and plane quintics, in the sense that we found an explicit polynomial formula for the Cayley-Salmon invariant of such curves. This formula is expressed as a determinant of a finite size matrix with polynomial entries; it is therefore very convenient for practical calculations and allows to determine in reasonable time and space whether the curve has undulation points or not.

The existence of such a formula rises several interesting questions:

\paragraph{$\bullet$} In this paper a "pedestrian" method of studying the structure of the undulation ideal was used. While it gives the correct answer, it does not give neither simple explanation of why such polynomials exist, nor deep insight into the structure of the 21 (resp. 36) polynomials that constitute the determinantal formula. It is interesting to find out if the results of this paper can be rederived using more involved methods of algebraic geometry.

\paragraph{$\bullet$} It would be interesting to find a generalization of these formulae to higher $r > 5$, or to make sure that such a generalization does not exist.

\paragraph{$\bullet$} Another related direction of research is application of the method used in this paper to other types of invariants, associated with various other types of decomposition of curves. The undulation condition is associated with the decomposition
$$
P = a_1 a_1 a_1 a_1 + b_1 c_3
$$
where letters denote different polynomials and indices show their degrees. Similarly, one can consider different other types of decompositions, in particular,
\begin{align*}
& P = a_1 a_1 a_1 a_1 + b_2 c_2 \\
& P = a_2 b_1 b_1 + c_2 d_1 d_1 \\
& P = a_1 b_1 c_1 d_1 + a_1 b_1 c_1 e_1 + a_1 b_1 d_1 e_1 + a_1 c_1 d_1 e_1 + b_1 c_1 d_1 e_1 \\
& P = a_1 a_1 a_1 a_1 + b_1 b_1 b_1 b_1 + c_1 c_1 c_1 c_1 + d_1 d_1 d_1 d_1 + e_1 e_1 e_1 e_1 \\
& P = a_1 a_1 a_1 a_1 + b_2 b_2 + c_2 c_2
\end{align*}
It is easy to show that existence of each of these decompositions is equivalent to vanishing of a certain invariant polynomial in coefficients of $P$. Some of these decompositions were discussed already for a long time; in particular, the third decomposition
$$
P = a_1 b_1 c_1 d_1 + a_1 b_1 c_1 e_1 + a_1 b_1 d_1 e_1 + a_1 c_1 d_1 e_1 + b_1 c_1 d_1 e_1
$$
defines the so-called Lueroth quartics, and the corresponding invariant is called the Lueroth invariant. In analogue with the Cayley-Salmon undulation invariant, it has high degree -- 54 -- and no known explicit formula. It would be interesting to see whether the method of the current paper can produce a determinantal formula for these invariants.

\section*{Acknowledgements}

We are indebted to B.Sturmfels for letting us know about the undulation problem. We are grateful to A.Morozov and all participants of ITEP mathematical physics group seminar for stimulating discussions. Work is partly supported by RFBR grants 12-01-$33071\underline{\hspace{1ex}}{\rm mol}\underline{\hspace{1ex}}{\rm a}\underline{\hspace{1ex}}{\rm ved}$, 10-02-00509, 12-02-92108-YaF-a 11-01-92612-KO 11-02-90453-UKR (A.P), RFBR grant 12-01-00525 (Sh.Sh), grant for support of scientific schools NSh-3349.2012.2 and government contract 8606.

\pagebreak

\section*{Appendix: explicit formulae for the 21 polynomials}

In this section, we denote
$$
S_{400} = a, S_{310} = 4b, S_{301} = 4c, S_{220} = 6d, S_{211} = 12e, S_{202} = 6f
$$
$$
S_{130} = 4g, S_{121} = 12h, S_{112} = 12i, S_{103}=4j, S_{040}=k, S_{031}=4l, S_{022}=6m, S_{013}=4n, S_{004} = o
$$

\paragraph{Representation (1).}The explicit formulae for the 3 polynomials $\alpha_1, \alpha_2, \alpha_3$ that transform under the fundamental representation (1) of $SL(3)$ are
$$
\alpha_1 = v_1 \phi, \ \ \ \alpha_2 = v_2 \phi, \ \ \ \alpha_3 = v_3 \phi
$$
where
{\fontsize{8pt}{0pt}{
\begin{center}
$\phi = (k o-4 l n+3 m^2) v_{1}^4+(-4 o g-12 m i+4 j l+12 n h) v_{1}^3 v_{2}+(-12 m h+12 l i+4 n g-4 j k) v_{1}^3 v_{3}+(-12 j h-12 n e+6 m f+12 i^2+6 o d) v_{1}^2 v_{2}^2+(-12 l f-12 h i+24 m e-12 n d+12 j g) v_{1}^2 v_{2} v_{3}+(-12 g i+6 k f-12 l e+6 m d+12 h^2) v_{1}^2 v_{3}^2+(4 n c-4 o b+12 j e-12 f i) v_{1} v_{2}^3+(-12 m c-12 e i-12 j d+12 n b+24 f h) v_{1} v_{2}^2 v_{3}+(-12 m b-12 e h+12 l c-12 f g+24 d i) v_{1} v_{2} v_{3}^2+(4 l b-12 d h+12 e g-4 k c) v_{1} v_{3}^3+(a o+3 f^2-4 c j) v_{2}^4+(12 i c-12 f e-4 a n+4 b j) v_{2}^3 v_{3}+(6 a m-12 i b+12 e^2+6 f d-12 h c) v_{2}^2 v_{3}^2+(12 h b+4 g c-12 e d-4 a l) v_{2} v_{3}^3+(-4 g b+3 d^2+a k) v_{3}^4$
\end{center}}}
\hspace{-4ex} is the unique basis polynomial of ${\cal I}_{2,4}$.

\paragraph{Representation (1,1).}The explicit formulae for the 3 polynomials $\beta_1, \beta_2, \beta_3$ that transform under the antisymmetric representation (1,1) of $SL(3)$ are

{\fontsize{8pt}{0pt}{
\begin{center}
$\beta_1 = (-3 l^2 o+6 l m n+3 k m o-3 k n^2-3 m^3) v_{1}^5+(-10 g m o+5 k j n-15 h m n+15 i m^2-5 k i o+10 g n^2+15 h l o-10 l i n-5 l m j) v_{1}^4 v_{2}+(10 l^2 j+15 h m^2-5 k h o-15 l i m-5 g m n+5 g l o+15 k i n-10 h l n-10 k m j) v_{1}^4 v_{3}+(24 h i n-12 d n^2-6 f m^2+12 d m o-2 k j^2+12 e m n+4 l f n+16 g i o+8 l i j-12 e l o-16 g j n-24 i^2 m+2 k f o-18 h^2 o+12 h m j) v_{1}^3 v_{2}^2+(4 g h o+24 h^2 n+4 k i j-40 h l j-12 h i m-12 d l o+12 d m n+16 e l n+36 g m j-24 e m^2+12 l f m-40 g i n-12 k f n+8 k e o+24 l i^2) v_{1}^3 v_{2} v_{3}+(2 d k o-6 d m^2+12 e l m-18 k i^2+24 h l i+16 k h j+12 g i m+8 g h n-12 l^2 f+4 d l n-2 g^2 o+12 k f m-16 g l j-12 k e n-24 h^2 m) v_{1}^3 v_{3}^2+(-9 e m j+3 c l o-3 c m n-9 h f n+6 b n^2+18 d j n+6 g j^2-3 l f j-18 h i j-18 e i n+12 i^3-6 b m o-18 d i o+27 e h o+18 f i m-6 g f o) v_{1}^2 v_{2}^3+(-18 g i j-15 g e o-27 l f i+9 b l o+9 d h o+33 g f n+36 h^2 j-9 h f m+3 k f j-54 e h n-18 h i^2-9 b m n+45 e i m+24 e l j+36 d i n-3 k c o-45 d m j+9 c m^2-6 c l n) v_{1}^2 v_{2}^2 v_{3}+(36 g i^2+45 e h m-9 d i m-18 h^2 i+24 g e n+3 d g o+36 h l f-15 k e j-18 g h j+9 k f i+33 d l j-27 d h n-54 e l i-45 g f m-6 b l n-3 b k o-9 c l m+9 b m^2+9 k c n) v_{1}^2 v_{2} v_{3}^2+(12 h^3+6 g^2 j-6 d k j+18 g l f-3 b l m+3 b k n-18 k h f-6 k c m-18 g h i+6 c l^2-18 e h l+18 d h m-3 d g n-9 d l i-9 g e m+27 k e i) v_{1}^2 v_{3}^3+(-9 e^2 o+12 e i j+2 c m j-6 c h o+8 b i o+6 d f o-a n^2-3 f^2 m-6 d j^2+4 c i n+a m o+6 e f n-8 b j n+6 h f j-12 f i^2) v_{1} v_{2}^4+(-36 e h j+2 a m n+30 h f i-2 a l o-6 g f j-10 b h o+24 e^2 n+6 d e o+6 l f^2-4 c l j-30 d f n+16 c h n-12 e i^2+6 g c o+22 b m j+24 d i j-12 b i n-6 e f m-18 c i m) v_{1} v_{2}^3 v_{3}+(-20 b l j+24 g e j-3 d^2 o+48 d f m+2 k c j-6 e l f-6 g f i+2 b g o-6 d e n-3 a m^2+a k o+24 b h n+24 c l i-30 d i^2-6 c h m-6 b i m-6 d h j-30 h^2 f-3 k f^2+2 a l n+48 e h i-30 e^2 m-20 g c n) v_{1} v_{2}^2 v_{3}^2+(-12 c h l+24 g h f-12 e h^2-4 b g n-2 a k n-18 b h m+6 d^2 n+16 b l i-6 d e m+30 d h i+2 a l m+6 k e f-6 d g j-30 d l f+22 g c m+24 e^2 l+6 b k j-36 g e i-10 k c i) v_{1} v_{2} v_{3}^3+(-9 k e^2+6 d g i+8 k c h-6 b k i+4 b h l+a k m+12 g e h-6 g^2 f-8 g c l+6 d e l-3 d^2 m+6 d k f-a l^2-12 d h^2+2 b g m) v_{1} v_{3}^4+(-2 b f o+3 f^2 i-2 c i j+3 c e o-a i o-3 e f j+2 b j^2+a j n-c f n) v_{2}^5+(2 a h o-6 h f^2-3 d c o+9 b f n+a i n+3 c f m+4 c h j+b e o-10 b i j+6 c i^2-10 c e n-3 e f i+3 d f j+6 e^2 j-3 a m j) v_{2}^4 v_{3}+(b d o-3 c l f-6 e^2 i+8 b h j+15 e h f-a g o-5 a h n+3 a i m-3 d f i-15 b f m-18 c h i+11 d c n+3 a l j+12 b i^2-2 g c j-6 b e n+3 g f^2+12 c e m-9 d e j) v_{2}^3 v_{3}^2+(-3 d h f+12 b e m-3 b d n-2 b g j-6 e^2 h+k c f-5 a l i+3 d^2 j+3 a h m-6 c e l+15 d e i-9 g e f-18 b h i+3 a g n-15 d c m+8 g c i+11 b l f+12 c h^2-a k j) v_{2}^2 v_{3}^3+(-6 d^2 i-3 a g m+a h l-3 d e h-3 b k f+6 b h^2+3 d g f+9 d c l+k c e+6 g e^2-10 g c h+4 b g i+2 a k i-10 b e l+3 b d m) v_{2} v_{3}^4+(-2 d k c-3 d g e+3 d^2 h-b d l+3 b k e+a g l-2 b g h-a k h+2 g^2 c) v_{3}^5$
\end{center}
}}

{\fontsize{8pt}{0pt}{
\begin{center}
$\beta_2 =
(3 h m n-2 g n^2-3 i m^2+l m j+2 l i n-k j n+2 g m o+k i o-3 h l o) v_{1}^5+(12 i^2 m-6 h m j-12 h i n+6 e l o-k f o+9 h^2 o+8 g j n+6 d n^2-2 l f n+k j^2-6 e m n-8 g i o+3 f m^2-6 d m o-4 l i j) v_{1}^4 v_{2}+(-4 e l n-2 k e o+6 e m^2+3 d l o+10 g i n+10 h l j-3 d m n-3 l f m-9 g m j-k i j-6 l i^2+3 h i m+3 k f n-6 h^2 n-g h o) v_{1}^4 v_{3}+(-27 e h o-6 g j^2-3 c l o+18 e i n+9 e m j-6 b n^2+3 l f j-18 f i m+6 g f o-18 d j n-12 i^3+18 d i o+6 b m o+9 h f n+3 c m n+18 h i j) v_{1}^3 v_{2}^2+(36 e h n-6 d h o+10 g e o+12 h i^2+6 h f m+2 k c o+30 d m j+6 b m n-24 d i n+4 c l n-16 e l j+12 g i j-6 b l o-30 e i m-2 k f j+18 l f i-6 c m^2-24 h^2 j-22 g f n) v_{1}^3 v_{2} v_{3}+(b k o-3 k f i+9 d h n-12 g i^2-15 e h m+6 h^2 i+2 b l n+18 e l i+3 c l m-d g o-12 h l f+6 g h j-3 k c n-11 d l j+3 d i m+15 g f m-8 g e n+5 k e j-3 b m^2) v_{1}^3 v_{3}^2+(-12 e f n+18 e^2 o+24 f i^2-4 c m j+6 f^2 m+2 a n^2-8 c i n+16 b j n+12 c h o-12 d f o-24 e i j-2 a m o-12 h f j+12 d j^2-16 b i o) v_{1}^2 v_{2}^3+(15 b h o-9 d e o-9 l f^2+6 c l j+54 e h j-36 e^2 n+18 b i n+3 a l o+27 c i m-33 b m j-36 d i j+9 g f j+9 e f m-3 a m n+18 e i^2+45 d f n-9 g c o-24 c h n-45 h f i) v_{1}^2 v_{2}^2 v_{3}+(3 k f^2+6 d e n-a k o-24 g e j+6 d h j+20 b l j-48 d f m+6 e l f-2 k c j+30 d i^2+6 c h m-2 a l n-2 b g o+6 b i m+3 d^2 o+20 g c n-24 c l i+6 g f i-24 b h n-48 e h i+3 a m^2+30 h^2 f+30 e^2 m) v_{1}^2 v_{2} v_{3}^2+(6 e h^2+2 b g n-a l m-3 k e f+9 b h m+18 g e i+15 d l f+a k n+6 c h l-3 d^2 n-3 b k j+3 d e m-15 d h i-12 g h f+3 d g j-8 b l i-12 e^2 l-11 g c m+5 k c i) v_{1}^2 v_{3}^3+(-10 b j^2-15 c e o-5 a j n+5 c f n+5 a i o+10 c i j-15 f^2 i+15 e f j+10 b f o) v_{1} v_{2}^4+(12 e f i+24 h f^2-12 d f j-16 c h j+40 b i j-4 a i n+12 d c o-36 b f n+12 a m j-4 b e o-12 c f m-24 c i^2-24 e^2 j+40 c e n-8 a h o) v_{1} v_{2}^3 v_{3}+(-33 d c n-24 b h j+9 d f i+6 g c j-36 b i^2+27 d e j+3 a g o-9 a l j+45 b f m+18 e^2 i-3 b d o-9 a i m-9 g f^2-45 e h f+54 c h i+18 b e n+15 a h n-36 c e m+9 c l f) v_{1} v_{2}^2 v_{3}^2+(10 a l i+30 d c m-16 g c i+18 g e f-22 b l f-24 b e m-6 d^2 j+2 a k j+4 b g j+6 b d n-6 a g n-30 d e i-24 c h^2-2 k c f+36 b h i+6 d h f+12 c e l+12 e^2 h-6 a h m) v_{1} v_{2} v_{3}^3+(3 b k f+3 a g m-3 d g f-k c e-4 b g i+10 b e l-6 g e^2-a h l-6 b h^2-3 b d m+6 d^2 i+3 d e h-2 a k i-9 d c l+10 g c h) v_{1} v_{3}^4+(-6 c f j+3 a j^2+3 c^2 o-3 a f o+3 f^3) v_{2}^5+(-15 e f^2-15 a i j-5 b c o+15 c f i+10 c e j-10 c^2 n+10 a f n+5 b f j+5 a e o) v_{2}^4 v_{3}+(-8 b e j-2 a d o+6 d f^2-12 b f i+24 e^2 f-4 d c j-12 a f m+2 b^2 o-16 a e n+18 a i^2-24 c e i-12 c h f+16 b c n+12 a h j+12 c^2 m) v_{2}^3 v_{3}^2+(6 a d n+9 b h f+3 g c f-12 e^3-6 c^2 l+18 a e m-6 b^2 n-27 a h i-18 d e f+9 d c i+18 b e i+18 c e h+3 b d j+6 a l f-18 b c m-3 a g j) v_{2}^2 v_{3}^3+(-6 d c h+6 a g i-a k f+12 d e^2-4 g c e+6 b^2 m-2 b g f+k c^2+9 a h^2-6 a d m+3 d^2 f-12 b e h+8 b c l-6 b d i-8 a e l) v_{2} v_{3}^4+(-2 b^2 l+a k e-3 a g h-b k c+2 a d l-3 d^2 e+3 b d h+d g c+2 b g e) v_{3}^5
$
\end{center}
}}

{\fontsize{8pt}{0pt}{
\begin{center}
$\beta_3 =
(-2 k m j+g l o+2 l^2 j-g m n-k h o+3 k i n-2 h l n+3 h m^2-3 l i m) v_{1}^5+(9 g m j-10 h l j+3 l f m+6 h^2 n+4 e l n-3 h i m+2 k e o-6 e m^2-3 d l o+6 l i^2+k i j+g h o-3 k f n+3 d m n-10 g i n) v_{1}^4 v_{2}+(-6 l^2 f+6 e l m+4 g h n+6 k f m-9 k i^2+12 h l i-8 g l j-12 h^2 m+8 k h j+d k o+6 g i m-3 d m^2+2 d l n-g^2 o-6 k e n) v_{1}^4 v_{3}+(12 h^2 j-2 c l n+3 d h o-5 g e o-6 h i^2+11 g f n-6 g i j+k f j-3 h f m+3 b l o-3 b m n-15 d m j-18 e h n-k c o+15 e i m+8 e l j+3 c m^2+12 d i n-9 l f i) v_{1}^3 v_{2}^2+(-6 d i m-10 k e j-18 d h n-12 h^2 i-4 b l n+24 h l f-2 b k o+6 b m^2-6 c l m+22 d l j+6 k f i-36 e l i+2 d g o+16 g e n+6 k c n+24 g i^2-12 g h j+30 e h m-30 g f m) v_{1}^3 v_{2} v_{3}+(12 h^3+6 g^2 j-6 d k j+18 g l f-3 b l m+3 b k n-18 k h f-6 k c m-18 g h i+6 c l^2-18 e h l+18 d h m-3 d g n-9 d l i-9 g e m+27 k e i) v_{1}^3 v_{3}^2+(-6 e i^2+12 d i j-3 e f m+15 h f i+3 g c o-2 c l j-6 b i n+3 l f^2-18 e h j-5 b h o-15 d f n+3 d e o+11 b m j+a m n+12 e^2 n-9 c i m+8 c h n-3 g f j-a l o) v_{1}^2 v_{2}^3+(-20 b l j+24 g e j-3 d^2 o+48 d f m+2 k c j-6 e l f-6 g f i+2 b g o-6 d e n-3 a m^2+a k o+24 b h n+24 c l i-30 d i^2-6 c h m-6 b i m-6 d h j-30 h^2 f-3 k f^2+2 a l n+48 e h i-30 e^2 m-20 g c n) v_{1}^2 v_{2}^2 v_{3}+(-9 d g j-45 d l f-3 a k n+9 k e f+45 d h i+3 a l m+33 g c m+36 e^2 l-27 b h m-18 c h l-6 b g n-54 g e i+36 g h f+24 b l i+9 d^2 n+9 b k j-15 k c i-9 d e m-18 e h^2) v_{1}^2 v_{2} v_{3}^2+(24 g e h+4 b g m-16 g c l-12 b k i-18 k e^2+12 d g i-24 d h^2-2 a l^2-12 g^2 f+2 a k m+12 d e l+16 k c h-6 d^2 m+8 b h l+12 d k f) v_{1}^2 v_{3}^3+(2 a h o-6 h f^2-3 d c o+9 b f n+a i n+3 c f m+4 c h j+b e o-10 b i j+6 c i^2-10 c e n-3 e f i+3 d f j+6 e^2 j-3 a m j) v_{1} v_{2}^4+(-4 g c j-18 d e j-6 c l f+24 c e m+6 a l j+6 g f^2-30 b f m-12 b e n+24 b i^2+6 a i m-36 c h i-10 a h n-2 a g o+30 e h f+2 b d o-12 e^2 i-6 d f i+22 d c n+16 b h j) v_{1} v_{2}^3 v_{3}+(-9 d h f+9 d^2 j+9 a h m-6 b g j-27 g e f+9 a g n+36 c h^2-18 e^2 h+36 b e m+45 d e i-54 b h i+24 g c i+33 b l f-45 d c m-3 a k j-15 a l i-9 b d n+3 k c f-18 c e l) v_{1} v_{2}^2 v_{3}^2+(8 a k i+12 d g f-12 a g m-12 b k f+4 k c e+12 b d m-24 d^2 i-12 d e h+24 g e^2+4 a h l-40 b e l+24 b h^2-40 g c h+16 b g i+36 d c l) v_{1} v_{2} v_{3}^3+(-15 d g e+15 b k e-10 d k c-5 a k h+10 g^2 c+5 a g l-5 b d l+15 d^2 h-10 b g h) v_{1} v_{3}^4+(2 c^2 n+3 e f^2-2 a f n+3 a i j-2 c e j-b f j+b c o-3 c f i-a e o) v_{2}^5+(-6 a h j+6 a f m-6 c^2 m+12 c e i+8 a e n-8 b c n-b^2 o-9 a i^2+4 b e j-12 e^2 f+a d o+6 c h f-3 d f^2+6 b f i+2 d c j) v_{2}^4 v_{3}+(-18 a e m+18 d e f+18 b c m+6 c^2 l+27 a h i-18 b e i-6 a d n-9 b h f-3 g c f+6 b^2 n-9 d c i-3 b d j+3 a g j-18 c e h+12 e^3-6 a l f) v_{2}^3 v_{3}^2+(-12 b^2 m+12 b d i+24 b e h-24 d e^2-18 a h^2+4 b g f-12 a g i+16 a e l-6 d^2 f-2 k c^2-16 b c l+8 g c e+12 a d m+2 a k f+12 d c h) v_{2}^2 v_{3}^3+(-15 b d h-10 b g e+10 b^2 l-5 a k e-5 d g c+15 a g h-10 a d l+5 b k c+15 d^2 e) v_{2} v_{3}^4+(3 a d k-3 d^3-3 a g^2+6 b d g-3 b^2 k) v_{3}^5
$
\end{center}
}}

\paragraph{Representation (3,2).} The explicit formulae for the 15 polynomials $\beta_4, \ldots, \beta_{18}$ that transform under the irreducible representation (3,2) of $SL(3)$ are

{\fontsize{8pt}{0pt}{
\begin{center}$
\beta_4 = (k m o-m^3-l^2 o+2 l m n-k n^2) v_{1}^4 v_{3}+(2 h m o+2 i m n-2 m^2 j+2 l j n-2 l i o-2 h n^2) v_{1}^3 v_{2}^2+(-4 h m n-4 l i n+4 h l o+4 i m^2+4 g n^2-4 g m o) v_{1}^3 v_{2} v_{3}+(-2 k h o+4 k i n+2 l^2 j+2 g l o-2 k m j+4 h m^2-2 g m n-2 h l n-4 l i m) v_{1}^3 v_{3}^2+(2 l f o-4 e m o+2 h i o-4 i^2 n-2 h j n-2 l j^2+6 i m j+4 e n^2-2 f m n) v_{1}^2 v_{2}^3+(-k f o+2 g i o-5 d n^2+4 f m^2+16 h i n+e l o-2 h m j+5 d m o-2 g j n+k j^2+2 l i j-10 i^2 m-3 l f n-e m n-7 h^2 o) v_{1}^2 v_{2}^2 v_{3}+(-2 h i m+4 e l n-12 g i n+4 g h o+5 d m n+8 g m j-l f m+8 l i^2-5 e m^2+k e o-5 d l o+k f n-6 h l j-2 k i j+2 h^2 n) v_{1}^2 v_{2} v_{3}^2+(-3 k e n+d k o-l^2 f+4 g i m-g^2 o-2 d m^2+k f m+3 e l m-4 h^2 m-4 g l j+2 h l i+2 g h n+4 k h j+d l n-3 k i^2) v_{1}^2 v_{3}^3+(2 e i o+4 h j^2-2 f m j-2 e j n-2 c n^2+2 c m o-4 h f o-4 i^2 j+6 f i n) v_{1} v_{2}^4+(-2 g j^2+7 e h o+3 h f n-10 h i j+12 i^3+2 b n^2-2 b m o-14 e i n+4 d j n-3 c l o-4 d i o+7 e m j+3 l f j+2 g f o-8 f i m+3 c m n) v_{1} v_{2}^3 v_{3}+(-6 e l j-3 b m n-k f j+d h o-8 e h n+3 b l o+l f i-3 c m^2+19 e i m+h f m+2 c l n-5 g e o+6 g i j-18 h i^2+10 d i n+k c o-11 d m j+12 h^2 j-g f n) v_{1} v_{2}^2 v_{3}^2+(-k c n-10 g h j-2 b l n+9 d l j+e h m+4 g i^2+8 g e n+3 b m^2-10 e l i-5 d h n+k f i-b k o+d g o-5 d i m+2 h l f+6 h^2 i-3 g f m+c l m+k e j) v_{1} v_{2} v_{3}^3+(-2 g h i+2 g^2 j-b l m-2 k h f+4 d h m-d g n+2 g l f-d l i+3 k e i-3 g e m+b k n-2 d k j) v_{1} v_{3}^4+(2 e f o-2 f^2 n-2 e j^2+2 c j n+2 f i j-2 c i o) v_{2}^5+(-5 c m j+2 b i o-3 h f j+d j^2-4 e^2 o+3 c h o-6 f i^2+2 c i n-2 b j n-d f o+5 f^2 m+8 e i j) v_{2}^4 v_{3}+(8 e^2 n-4 d i j+3 c i m+4 c l j+4 d e o-6 c h n-4 l f^2+15 h f i+g f j+5 b m j-2 b i n-8 e f m-g c o-6 e i^2-3 b h o-6 e h j) v_{2}^3 v_{3}^2+(-3 b i m+3 d i^2-4 c l i-d^2 o+b g o+k f^2+4 d f m-9 h^2 f-4 b l j+2 g e j+6 e h i-8 d e n+2 g c n-4 e^2 m+8 e l f-5 g f i+3 c h m-k c j+6 b h n+3 d h j) v_{2}^2 v_{3}^3+(-2 k e f-2 b g n-3 b h m-g c m+b k j-2 g e i-d g j+4 d e m-4 d l f+2 d^2 n+k c i+6 g h f-3 d h i+4 b l i) v_{2} v_{3}^4+(d g i+d k f-d^2 m+b g m-b k i-g^2 f) v_{3}^5$
\end{center}
}}

{\fontsize{8pt}{0pt}{
\begin{center}$
\beta_5 = (4 i m n+4 l j n-4 m^2 j+4 h m o-4 h n^2-4 l i o) v_{1}^4 v_{2}+(-2 i m^2+2 h m n-2 k j n-2 h l o+2 l m j+2 k i o) v_{1}^4 v_{3}+(4 h i o-4 l j^2+4 l f o+12 i m j-4 h j n-4 f m n+8 e n^2-8 i^2 n-8 e m o) v_{1}^3 v_{2}^2+(-4 i^2 m+10 e l o+2 d n^2+4 g j n+12 f m^2-10 l f n-8 h m j-2 d m o+12 h i n-4 g i o-10 e m n-2 k f o-2 h^2 o+2 k j^2) v_{1}^3 v_{2} v_{3}+(-5 l f m+3 e m^2+d l o+4 h l j-2 g m j+2 g h o-4 h^2 n-d m n+2 h i m+5 k f n-2 k i j-3 k e o) v_{1}^3 v_{3}^2+(-4 f m j+8 h j^2-8 h f o+4 c m o-4 e j n-4 c n^2+4 e i o+12 f i n-8 i^2 j) v_{1}^2 v_{2}^3+(16 i^3+8 l f j-12 e i n-4 b n^2+8 c m n+16 e m j-12 h i j-4 g j^2+16 h f n-4 d j n+4 g f o-8 c l o-4 e h o-28 f i m+4 b m o+4 d i o) v_{1}^2 v_{2}^2 v_{3}+(3 k c o-20 e l j-d h o+14 e i m-4 b l o-6 d i n+7 d m j+8 g i j+4 e h n+8 c l n-10 g f n+8 h^2 j-3 h f m+16 l f i+4 b m n+2 g e o-11 c m^2-3 k f j-16 h i^2) v_{1}^2 v_{2} v_{3}^2+(-2 d l j-b m^2+6 k e j-2 h l f-3 k f i+4 c l m-4 k c n-d g o+5 g f m+4 h^2 i+4 d h n-d i m+b k o-6 e h m-4 g h j) v_{1}^2 v_{3}^3+(-4 c i o+4 c j n+4 f i j+4 e f o-4 e j^2-4 f^2 n) v_{1} v_{2}^4+(-8 e f n-4 b i o+4 b j n+12 e i j-2 a m o+2 d j^2-2 e^2 o-10 c m j-2 d f o-4 f i^2+12 f^2 m+10 c h o+2 a n^2-10 h f j) v_{1} v_{2}^3 v_{3}+(8 c l j+4 e h j+3 a l o-4 g c o-10 b m j+4 g f j+16 c i m-6 d i j-20 c h n+8 b i n+7 d f n+14 h f i-d e o-16 e i^2-11 l f^2+2 b h o-3 a m n+8 e^2 n-3 e f m) v_{1} v_{2}^2 v_{3}^2+(-2 d h j-6 d f m+8 g c n-6 e^2 m-4 b h n+3 k f^2+8 e h i-2 a l n+3 a m^2-a k o-16 c l i+8 b l j+10 c h m+8 d i^2-4 g e j-6 h^2 f-4 g f i-2 d e n-4 b i m+10 e l f+d^2 o-2 k c j) v_{1} v_{2} v_{3}^3+(-4 g c m+2 d g j+4 k c i-a l m+a k n+2 g h f+2 b h m+d l f+3 d e m-2 b k j-4 d h i-3 k e f-d^2 n) v_{1} v_{3}^4+(2 a i o-2 f^2 i+2 c f n-2 a j n-2 c e o+2 e f j) v_{2}^4 v_{3}+(3 h f^2-5 c f m-4 e^2 j-3 a h o+5 a m j+2 e f i+2 b e o+4 c e n-d f j+d c o-2 b f n-2 a i n) v_{2}^3 v_{3}^2+(-2 c e m-3 a i m-d f i-6 e h f-4 b e n+a g o-b d o-4 a l j+4 d e j-2 d c n+5 b f m+4 c l f-g f^2+4 e^2 i+6 a h n) v_{2}^2 v_{3}^3+(a k j-3 a h m+2 g e f+2 b d n+2 b e m+4 a l i+d c m+3 d h f-d^2 j-k c f-4 d e i-4 b l f-2 a g n) v_{2} v_{3}^4+(d^2 i-a k i-b d m+b k f-d g f+a g m) v_{3}^5
$
\end{center}
}}

{\fontsize{8pt}{0pt}{
\begin{center}$
\beta_6 = (2 h m o+2 i m n-2 m^2 j+2 l j n-2 l i o-2 h n^2) v_{1}^5+(2 l f o-4 e m o+2 h i o-4 i^2 n-2 h j n-2 l j^2+6 i m j+4 e n^2-2 f m n) v_{1}^4 v_{2}+(3 e l o+2 l i j-3 e m n-5 l f n-2 g j n-4 h^2 o+2 g i o+d n^2-d m o+8 h i n+5 f m^2-6 i^2 m) v_{1}^4 v_{3}+(2 e i o+4 h j^2-2 f m j-2 e j n-2 c n^2+2 c m o-4 h f o-4 i^2 j+6 f i n) v_{1}^3 v_{2}^2+(2 b m o-10 e i n-4 d i o+4 d j n-8 f i m+3 l f j-14 h i j+7 e h o-2 b n^2+12 i^3+7 h f n+3 e m j-2 g f o+3 c m n+2 g j^2-3 c l o) v_{1}^3 v_{2} v_{3}+(-6 h i^2+b m n+15 e i m-4 d i n-2 g i j+4 d h o-6 e h n-8 h f m+4 c l n-6 e l j-b l o+3 l f i-3 g e o+8 h^2 j+5 g f n-4 c m^2) v_{1}^3 v_{3}^2+(2 e f o-2 f^2 n-2 e j^2+2 c j n+2 f i j-2 c i o) v_{1}^2 v_{2}^3+(-7 e^2 o-3 c m j-5 d j^2-h f j+c h o+2 c i n-2 e f n-2 b j n+a n^2+4 f^2 m-a m o-10 f i^2+5 d f o+2 b i o+16 e i j) v_{1}^2 v_{2}^2 v_{3}+(3 g c o-6 c h n+a l o+2 c l j+d e o+10 d i j+c i m-3 g f j-a m n-3 l f^2-5 b h o+e f m-18 e i^2-b m j+12 e^2 n+6 b i n-11 d f n-8 e h j+19 h f i) v_{1}^2 v_{2} v_{3}^2+(2 b l j-a l n-4 h^2 f+3 e l f-4 c l i-5 b i m-9 e^2 m+a m^2+3 d e n+2 b h n+4 d f m+3 d i^2-d^2 o-8 d h j+b g o+8 c h m+6 g e j-4 g c n-3 g f i+6 e h i) v_{1}^2 v_{3}^3+(4 b j^2-4 b f o-4 e f j+4 f^2 i+4 c e o-4 c i j) v_{1} v_{2}^3 v_{3}+(a h o-2 e f i-6 c e n+2 e^2 j-5 d c o-12 b i j+8 c i^2+8 b f n+4 c h j-c f m-2 a i n+5 d f j-5 h f^2+4 b e o+a m j) v_{1} v_{2}^2 v_{3}^2+(-10 c h i-a g o+3 g f^2-3 b f m-5 d f i-2 g c j+8 b h j+6 e^2 i+a h n+e h f-10 b e n-5 d e j+2 c e m+a i m-a l j+9 d c n+b d o+4 b i^2+c l f) v_{1} v_{2} v_{3}^3+(-3 d e i-2 a h m-4 d c m+a l i+6 b e m-2 b h i-b d n+4 g c i-2 b g j-3 g e f+4 d h f+2 d^2 j+a g n-b l f) v_{1} v_{3}^4+(2 c f j-f^3-a j^2+a f o-c^2 o) v_{2}^4 v_{3}+(-2 c e j+2 c^2 n-2 a f n+4 a i j+2 b c o-2 a e o-4 c f i+4 e f^2-2 b f j) v_{2}^3 v_{3}^2+(-4 b c n-4 e^2 f+3 c h f+a f m-2 d f^2-3 a h j-3 a i^2+a d o-b^2 o-c^2 m+2 c e i+4 b f i+4 a e n+d c j+2 b e j) v_{2}^2 v_{3}^3+(4 d e f-2 a d n-d c i+3 a h i+a g j-3 b h f+2 b^2 n-2 b e i+2 b c m-2 a e m-g c f-b d j) v_{2} v_{3}^4+(-d^2 f+b d i+b g f-b^2 m-a g i+a d m) v_{3}^5$
\end{center}
}}

{\fontsize{8pt}{0pt}{
\begin{center}$
\beta_7 = (-2 g n^2+2 g m o+2 l i n+k j n-h l o-l m j+h m n-k i o-i m^2) v_{1}^4 v_{2}+(k i n-l i m+2 l^2 j+k h o-2 k m j+h m^2-g l o+g m n-2 h l n) v_{1}^4 v_{3}+(3 h^2 o-l f n+4 i^2 m+e l o-8 h i n-5 d m o+2 h m j-k j^2+k f o-e m n+5 d n^2) v_{1}^3 v_{2}^2+(4 l f m-4 g h o+4 h^2 n-4 l i^2-4 d m n+4 k i j+4 g i n-4 k f n+4 d l o-4 h l j) v_{1}^3 v_{2} v_{3}+(-3 k i^2+5 k f m+d l n-4 h^2 m+e l m-5 l^2 f-2 g i m-d k o-k e n+g^2 o+8 h l i) v_{1}^3 v_{3}^2+(-4 i^3-4 b n^2-4 d j n+4 h f n-2 f i m-2 e m j+4 d i o+2 h i j+2 g j^2+4 b m o+8 e i n-2 g f o-6 e h o) v_{1}^2 v_{2}^3+(9 d m j-k c o+4 b m n+4 g e o+d h o+14 h i^2-4 h^2 j-4 b l o-7 h f m-10 d i n+c m^2-10 g i j+k f j-8 e i m+6 g f n+4 e l j) v_{1}^2 v_{2}^2 v_{3}+(-4 c l m+7 d i m-4 g e n+b k o+4 g i^2-b m^2+10 g h j-9 g f m+8 e h m-6 d l j-d g o+4 k c n-4 k e j+10 h l f-14 h^2 i-k f i) v_{1}^2 v_{2} v_{3}^2+(-8 e h l-4 k c m+4 h^3-4 k h f-2 g^2 j+4 c l^2+4 g l f+2 d h m-4 d l i+6 k e i+2 g e m+2 d k j-2 g h i) v_{1}^2 v_{3}^3+(-a m o+d f o-d j^2-h f j+c h o-2 e i j+4 f i^2+3 e^2 o+4 b j n+a n^2-4 e f n+c m j-4 b i o-2 c i n) v_{1} v_{2}^4+(9 e f m-g f j-5 d e o-2 c l j-a m n-9 b m j+c i m+6 e h j-2 e i^2-d f n+6 d i j+l f^2-9 h f i+6 b i n-4 e^2 n+a l o+g c o+3 b h o) v_{1} v_{2}^3 v_{3}+(-b i m+7 g f i+4 c l i+c h m-4 b h n-10 e l f-k f^2-5 d i^2-7 d h j-6 g c n+d^2 o+6 b l j+10 d e n+k c j-b g o+5 h^2 f) v_{1} v_{2}^2 v_{3}^2+(-9 d e m-6 c h l-b k j+9 g c m+a l m+9 d h i+d g j+5 k e f-6 g h f-6 g e i-d^2 n+4 e^2 l+2 b g n-3 k c i-b h m-a k n+2 e h^2+d l f) v_{1} v_{2} v_{3}^3+(4 k c h-4 d h^2+d g i+g^2 f+2 g e h-a l^2+a k m+2 b h l-b k i-b g m-4 g c l-d k f-3 k e^2+4 d e l) v_{1} v_{3}^4+(c f n-a j n-f^2 i+e f j-c e o+a i o) v_{2}^5+(-2 e^2 j+h f^2-2 c f m-a i n-a h o+2 a m j+2 c e n+b e o-b f n+e f i) v_{2}^4 v_{3}+(-2 b e n-e h f-a l j-2 c e m+2 e^2 i-2 d f i+d c n+a h n+2 b f m+c l f+d e j) v_{2}^3 v_{3}^2+(a g n-2 d c m+2 d h f+d e i-b d n-2 e^2 h-b l f-g e f+2 b e m-a l i+2 c e l) v_{2}^2 v_{3}^3+(-2 a g m-d^2 i+d c l+a k i+2 b d m+2 g e^2+a h l-k c e-d e h-2 b e l) v_{2} v_{3}^4+(-a k h+a g l+b k e-b d l-d g e+d^2 h) v_{3}^5$
\end{center}
}}

{\fontsize{8pt}{0pt}{
\begin{center}$
\beta_8 = (k m o-m^3-k n^2-l^2 o+2 l m n) v_{1}^5+(-2 h m n+2 h l o+2 i m^2-2 k i o+2 k j n-2 l m j) v_{1}^4 v_{2}+(4 k i n+4 h m^2-4 k m j-4 h l n-4 l i m+4 l^2 j) v_{1}^4 v_{3}+(k f o+l f n-2 f m^2-k j^2+3 d n^2+e m n-8 h i n-3 d m o-e l o+2 h m j+3 h^2 o+4 i^2 m) v_{1}^3 v_{2}^2+(-6 k f n-4 g h o-2 e m^2-4 l i^2+4 k i j-4 h l j+2 k e o+4 g i n-2 d m n+6 l f m+4 h^2 n+2 d l o) v_{1}^3 v_{2} v_{3}+(3 e l m-3 k i^2-7 l^2 f+3 d l n-d k o-4 h^2 m+g^2 o-2 g i m-2 d m^2+7 k f m+8 h l i-3 k e n) v_{1}^3 v_{3}^2+(-8 i^3-e m j-l f j+10 e i n+h f n-6 d j n+8 h i j+2 b m o+6 d i o-9 e h o-c m n+c l o-2 b n^2) v_{1}^2 v_{2}^3+(-2 k c o-4 g i j-2 h f m-14 d i n+14 e h n+4 c m^2+g f n+2 k f j+14 d m j-2 c l n+20 h i^2-16 h^2 j-b l o-21 e i m+4 e l j+b m n-l f i+3 g e o) v_{1}^2 v_{2}^2 v_{3}+(-11 d l j-4 g e n+16 g h j+6 h l f-5 k e j-5 e h m-8 g i^2+7 k c n-2 k f i-2 b l n-8 h^2 i+14 e l i+2 b m^2+12 d i m-d h n-4 g f m-7 c l m) v_{1}^2 v_{2} v_{3}^2+(-d g n-4 g^2 j-6 e h l+2 g l f+4 d k j+3 k e i-7 d l i+4 g h i+b k n+4 d h m-6 k c m+3 g e m-b l m+6 c l^2-2 k h f) v_{1}^2 v_{3}^3+(3 d j^2+c m j+4 b j n-f^2 m-3 d f o-h f j-2 c i n+8 f i^2-10 e i j-2 e f n-4 b i o+6 e^2 o+c h o) v_{1} v_{2}^4+(-3 d e o-20 e^2 n+6 b i n+7 e f m-g f j+14 e i^2-a l o+22 e h j+3 l f^2-2 c l j+13 d f n-25 h f i+a m n-10 d i j-9 b m j+g c o+3 b h o+c i m) v_{1} v_{2}^3 v_{3}+(4 c l i-2 k f^2+21 h^2 f-6 g c n+8 d e n+15 g f i-3 a m^2-8 g e j+c h m+a k o-b g o-40 e h i+11 d i^2+2 a l n-b i m+6 b l j-20 d f m-12 e l f+k c j+d h j-4 b h n+26 e^2 m) v_{1} v_{2}^2 v_{3}^2+(d^2 n+d g j-3 k c i-6 c h l+3 a l m-b k j-22 g h f+10 g e i+7 k e f+15 d l f+2 b g n-12 e^2 l-7 d h i+18 e h^2-11 d e m+9 g c m-b h m-3 a k n) v_{1} v_{2} v_{3}^3+(-2 a l^2+6 d e l+2 b h l-b k i-6 g e h+d g i-5 d k f-d^2 m+2 a k m+4 k c h-4 g c l+5 g^2 f-b g m) v_{1} v_{3}^4+(-2 b j^2+2 b f o+2 e f j-2 f^2 i+2 c i j-2 c e o) v_{2}^5+(d c o-6 c i^2-4 e^2 j-8 b f n+8 c e n-2 b e o-a m j-d f j+3 h f^2+10 b i j+2 e f i+c f m+a h o-4 c h j) v_{2}^4 v_{3}+(b d o+3 a i m-g f^2+4 e^2 i+2 g c j+18 c h i-d f i-14 c e m-4 a h n-a g o-6 e h f+4 d e j+2 a l j-12 b i^2+11 b f m-4 d c n-8 b h j+8 b e n-2 c l f) v_{2}^3 v_{3}^2+(-10 b e m-d^2 j+3 d h f-4 b d n+12 c e l+3 a h m+4 a g n+7 d c m-8 g c i+18 b h i-12 c h^2-6 a l i+k c f+2 g e f-6 b l f-a k j+2 b g j-4 d e i) v_{2}^2 v_{3}^3+(d^2 i-d g f+5 b d m+3 a k i-6 d c l+b k f+2 a h l-6 b h^2-4 k c e+4 b e l+10 g c h-5 a g m-4 b g i) v_{2} v_{3}^4+(-2 b d l+2 a g l+2 b g h-2 g^2 c-2 a k h+2 d k c) v_{3}^5$
\end{center}
}}

{\fontsize{8pt}{0pt}{
\begin{center}$
\beta_9 =(h l o+i m^2+k j n-h m n-k i o-l m j) v_{1}^5+(-k j^2+2 l i j-e l o+k f o-4 i^2 m-4 g j n+e m n+4 g i o-3 h^2 o-d m o-l f n+2 h i n+d n^2+4 h m j) v_{1}^4 v_{2}+(-2 h l j-2 k f n+k i j-h i m+2 l f m+g m j-e m^2+2 h^2 n-g h o+k e o) v_{1}^4 v_{3}+(-2 b n^2-8 h i j+4 i^3-4 d i o+2 b m o-2 e i n+2 h f n-4 e m j-4 g f o+4 g j^2+4 d j n+2 f i m+6 e h o) v_{1}^3 v_{2}^2+(-6 d i n+5 d h o+d m j-k c o+2 h i^2-9 h f m-6 e h n-l f i-c m^2-6 g i j+4 h^2 j+b m n+k f j-3 g e o-b l o+9 e i m+2 c l n+9 g f n) v_{1}^3 v_{2} v_{3}+(k c n-d h n+2 h l f+2 d i m-d l j-c l m+e h m-2 h^2 i-2 g f m-k e j+2 g h j) v_{1}^3 v_{3}^2+(h f j-5 d j^2-4 f i^2+c m j+a n^2-2 e f n-3 e^2 o+8 e i j-c h o-a m o+5 d f o) v_{1}^2 v_{2}^3+(-l f^2-14 e i^2-d e o+4 e^2 n-4 b h o-4 g f j-a m n+4 g c o+a l o-4 c h n+10 b i n-6 b m j-9 d f n+8 h f i+7 e f m+10 d i j) v_{1}^2 v_{2}^2 v_{3}+(-5 e^2 m-7 b i m+g f i-6 g c n+b g o+4 g e j-a l n+5 d i^2-e l f+10 c h m-10 d h j+6 b l j-d^2 o+7 d e n+a m^2-4 c l i) v_{1}^2 v_{2} v_{3}^2+(b h m+2 e h^2-2 d e m+k c i-2 g h f+d g j-d h i+g c m+2 d l f-b k j-2 c h l) v_{1}^2 v_{3}^3+(-a j n-2 c i j+a i o-e f j-2 b f o+2 b j^2+f^2 i+c f n+c e o) v_{1} v_{2}^4+(4 b e o-4 b i j-4 d c o+4 d f j-4 e^2 j-4 c f m-4 a i n+4 c e n+4 a m j+4 c i^2) v_{1} v_{2}^3 v_{3}+(g f^2+6 d c n+4 c l f-4 b i^2+b d o-10 b e n-7 d f i+4 b h j-10 c e m-4 a l j+a i m+14 e^2 i-a g o+9 b f m+4 a h n-8 e h f) v_{1} v_{2}^2 v_{3}^2+(3 a l i-2 e^2 h-4 c h^2-k c f+6 b h i-2 b g j+d^2 j+9 d h f-5 a h m-9 b l f-b d n+6 b e m-9 d e i+6 c e l+g e f-d c m+a g n+a k j) v_{1} v_{2} v_{3}^3+(d e h+a h l-2 d g f-a k i+2 g c h+d^2 i-k c e-2 b h^2+2 b k f-d c l) v_{1} v_{3}^4+(b c o-2 c^2 n+2 a f n-a e o-e f^2-a i j+c f i+2 c e j-b f j) v_{2}^4 v_{3}+(a h j+4 e^2 f+3 a i^2+2 b f i-c h f+a d o-d c j-8 c e i-5 a f m+5 c^2 m-b^2 o) v_{2}^3 v_{3}^2+(-4 b c m+4 d c i-2 d e f+2 b^2 n-6 a h i+2 b e i+4 a e m-2 b h f-4 c^2 l-2 a d n-4 e^3+4 a l f+8 c e h) v_{2}^2 v_{3}^3+(-2 g c e+4 b c l-a k f-4 a e l+a g i+3 a h^2-b d i-4 d c h-b^2 m-2 b e h+4 d e^2+k c^2+b g f+a d m) v_{2} v_{3}^4+(d g c-b k c-a g h+a k e-d^2 e+b d h) v_{3}^5
$
\end{center}
}}

{\fontsize{8pt}{0pt}{
\begin{center}$
\beta_{10} = (2 k j n-2 k i o+2 g m o-2 g n^2-2 l m j+2 l i n) v_{1}^5+(5 d n^2-6 h i n+2 k f o+4 g i o-e l o-f m^2+6 h m j-2 k j^2+e m n-4 g j n-l f n-5 d m o+2 l i j) v_{1}^4 v_{2}+(d l o-4 g h o+3 k e o+5 l f m-6 g m j+4 h l j+10 g i n-6 l i^2+e m^2+2 k i j-4 e l n-d m n-5 k f n) v_{1}^4 v_{3}+(-4 b n^2-7 e m j-l f j-c m n+2 d j n+4 b m o-2 d i o+3 e h o-6 h i j+c l o+4 f i m+4 e i n-6 g f o+3 h f n+6 g j^2) v_{1}^3 v_{2}^2+(15 d m j+b m n+9 g f n+7 d h o-12 h^2 j+3 k f j+18 h i^2+10 e h n-6 g i j-3 k c o-11 h f m-22 d i n+c m^2+2 c l n-l f i-3 g e o-7 e i m-b l o) v_{1}^3 v_{2} v_{3}+(-8 g e n-4 e h m-6 d l j-4 c l m-12 g i^2+d g o-10 h l f+7 g f m+2 b l n-6 k e j+3 k f i-b m^2+2 d h n-b k o+12 g h j+18 e l i+3 d i m+4 k c n) v_{1}^3 v_{3}^2+(a n^2-a m o-2 f^2 m+8 e i j+3 c m j-4 f i^2+7 d f o-3 e^2 o-3 c h o+3 h f j-7 d j^2-2 e f n) v_{1}^2 v_{2}^3+(-8 e^2 n-7 g f j+14 e h j-4 c h n-2 d e o-4 d f n-8 e i^2-2 c l j+7 g c o-5 h f i+12 e f m-c i m-5 b h o+16 b i n+2 l f^2+6 d i j-11 b m j) v_{1}^2 v_{2}^2 v_{3}+(-3 k f^2-8 b h n-a l n-4 c l i+6 b l j-2 d^2 o+15 d e n+21 d i^2+g f i+8 c h m+2 k c j+b i m-12 d h j-20 d f m-6 g c n+26 h^2 f+a k o-e l f+4 g e j+b g o-40 e h i+11 e^2 m) v_{1}^2 v_{2} v_{3}^2+(-2 d g j+11 d l f+4 e h^2-4 k c i-6 d h i+a l m-a k n-d e m-14 g h f+2 b k j-8 b l i+4 b h m-4 g c m+3 k e f+18 g e i-d^2 n+8 c h l+2 b g n-12 e^2 l) v_{1}^2 v_{3}^3+(-4 c i j+4 b j^2-4 e f j-4 b f o+4 c e o+4 f^2 i) v_{1} v_{2}^4+(4 c e n+2 a m j+6 d f j-2 h f^2-6 d c o-4 a i n-2 c f m+2 a h o-4 e^2 j+4 b e o+4 c i^2-4 b i j) v_{1} v_{2}^3 v_{3}+(4 g f^2-a l j-14 c e m+2 b d o+14 c h i-2 g c j+4 b h j+c l f-2 d f i+20 e^2 i+d c n+14 b f m-21 e h f-2 a g o+3 a h n-d e j-16 b i^2-4 b e n) v_{1} v_{2}^2 v_{3}^2+(g e f+k c f-20 c h^2+14 e^2 h+22 b h i-9 b l f+3 a l i+7 d h f+a g n-2 b g j+3 d^2 j+13 d c m+6 c e l-b d n-10 b e m-3 a h m-a k j-25 d e i) v_{1} v_{2} v_{3}^3+(2 d e h+a k i+3 d^2 i-4 b h^2-8 d c l-b k f+10 b e l-2 a h l-b d m+a g m-4 b g i-6 g e^2+d g f+8 g c h) v_{1} v_{3}^4+(a f o+2 c f j-f^3-c^2 o-a j^2) v_{2}^5+(-2 b f j+2 a i j+2 b c o-2 a e o-2 c f i+2 e f^2) v_{2}^4 v_{3}+(c h f-b^2 o-3 a f m+4 e^2 f-8 c e i+d c j+2 b f i+a d o+3 c^2 m-2 d f^2-a h j+3 a i^2) v_{2}^3 v_{3}^2+(d c i-2 c^2 l+6 a e m-b d j-8 e^3+8 b e i+2 a l f-b h f+10 c e h+a g j-g c f-9 a h i-6 b c m) v_{2}^2 v_{3}^3+(-10 b e h-b d i+6 a h^2+4 b c l-d^2 f-2 g c e-4 a e l+3 b^2 m-2 d c h+a g i+8 d e^2-3 a d m+b g f) v_{2} v_{3}^4+(2 b d h-2 a g h-2 d^2 e-2 b^2 l+2 a d l+2 b g e) v_{3}^5$
\end{center}
}}

{\fontsize{8pt}{0pt}{
\begin{center}$
\beta_{11} = (-e l o+l f n-d n^2+e m n-f m^2+d m o) v_{1}^5+(4 f i m-2 b m o-l f j+2 b n^2-3 h f n-2 d i o+2 d j n-c m n+c l o-2 e i n-e m j+3 e h o) v_{1}^4 v_{2}+(2 c m^2-2 c l n+6 d i n+b l o+g e o-2 e h n+4 e l j-4 d m j+4 h f m-2 d h o-3 l f i-g f n-3 e i m-b m n) v_{1}^4 v_{3}+(a m o+d f o-3 c h o+3 h f j-d j^2+4 b i o-a n^2+2 c i n-2 f^2 m-3 e^2 o+4 e f n-4 b j n+2 e i j+c m j-4 f i^2) v_{1}^3 v_{2}^2+(2 d i j+a m n+h f i-10 b i n+6 e i^2-5 e f m+b h o+d e o-10 e h j-g c o+4 e^2 n+8 c h n+9 b m j+3 l f^2+g f j-a l o-3 d f n-5 c i m-2 c l j) v_{1}^3 v_{2} v_{3}+(2 b h n-3 e l f-4 h^2 f-5 d e n-8 c h m+3 b i m-9 d i^2+4 d f m+3 g f i+6 c l i-4 b l j-a m^2+d^2 o+8 d h j+3 e^2 m+6 e h i-4 g e j-b g o+2 g c n+a l n) v_{1}^3 v_{3}^2+(4 c e o+2 b j^2-2 c f n-2 c i j-4 e f j+4 f^2 i-2 b f o-2 a i o+2 a j n) v_{1}^2 v_{2}^3+(-2 e f i-d f j-6 b i j+8 e^2 j-12 c e n-5 a m j+d c o+4 c h j+5 c f m+a h o-2 b e o-5 h f^2+8 b f n+4 a i n+2 c i^2) v_{1}^2 v_{2}^2 v_{3}+(2 g c j+12 b i^2-18 e^2 i-6 b h j-3 g f^2+10 c e m+3 a l j+19 e h f-b d o-d c n+6 b e n+d f i-11 b f m-3 c l f-5 a h n+a g o+a i m-8 c h i+d e j) v_{1}^2 v_{2} v_{3}^2+(5 b l f+4 a h m-4 d^2 j-8 d h f-3 a l i+b d n-6 g c i-4 b e m-6 e^2 h-a g n+4 b g j+15 d e i+8 c h^2-6 b h i+3 g e f-2 c e l) v_{1}^2 v_{3}^3+(a f o+2 c f j-f^3-c^2 o-a j^2) v_{1} v_{2}^4+(4 c^2 n+4 e f^2+4 a i j-4 a f n-4 c e j-4 c f i) v_{1} v_{2}^3 v_{3}+(-2 b f i-c h f-5 c^2 m-3 d c j+2 b e j-10 e^2 f+4 d f^2+2 a e n-2 b c n+5 a f m+16 c e i-7 a i^2+b^2 o-a d o+a h j) v_{1} v_{2}^2 v_{3}^2+(-2 b^2 n+2 c^2 l-3 a g j+3 b d j-4 a e m-2 a l f+12 e^3-14 c e h-8 d e f+2 a d n+3 g c f+4 b c m+3 d c i+7 b h f+7 a h i-10 b e i) v_{1} v_{2} v_{3}^3+(3 a g i+8 b e h-4 a h^2+5 d^2 f-3 b d i+2 g c e-a d m-2 b c l-5 b g f+2 a e l+b^2 m-6 d e^2) v_{1} v_{3}^4+(2 b c j+2 a f i-2 b f^2-2 c^2 i+2 c f e-2 a e j) v_{2}^3 v_{3}^2+(2 a d j-4 a f h-2 b^2 j+2 a e i+4 c^2 h-2 c d f-2 b c i+6 b e f-4 e^2 c) v_{2}^2 v_{3}^3+(-2 b c h-2 b d f+6 d c e+2 a g f-4 b e^2-2 g c^2+2 a e h+4 b^2 i-4 a d i) v_{2} v_{3}^4+(2 b g c-2 d^2 c+2 b d e-2 a g e-2 b^2 h+2 a d h) v_{3}^5
$
\end{center}
}}

{\fontsize{8pt}{0pt}{
\begin{center}$
\beta_{12} = (k m o-m^3-k n^2-l^2 o+2 l m n) v_{1}^4 v_{2}+(-2 g m o+4 h l o-4 h m n-2 l m j+2 g n^2+4 i m^2-2 l i n+2 k j n-2 k i o) v_{1}^3 v_{2}^2+(4 k i n+4 h m^2-4 k m j-4 h l n-4 l i m+4 l^2 j) v_{1}^3 v_{2} v_{3}+(-2 l^2 i+2 g l n+2 k i m-2 g m^2+2 h l m-2 k h n) v_{1}^3 v_{3}^2+(2 h i n-2 f m^2+3 e m n+2 l i j+4 h m j-4 i^2 m+l f n-3 h^2 o-k j^2+4 g i o-4 g j n+d m o-d n^2+k f o-3 e l o) v_{1}^2 v_{2}^3+(4 e l n-5 e m^2+5 l f m-6 g i n+d l o+4 k i j-d m n+8 g m j-5 k f n-2 g h o-2 h i m-12 h l j+k e o+8 h^2 n+2 l i^2) v_{1}^2 v_{2}^2 v_{3}+(5 k f m-10 h^2 m-7 k i^2+4 d m^2+2 k h j-5 l^2 f+k e n-2 g l j-e l m-2 g i m-3 d l n-d k o+2 g h n+16 h l i+g^2 o) v_{1}^2 v_{2} v_{3}^2+(-2 d l m-2 g l i+6 g h m-4 h^2 l-4 k e m+2 k h i+2 d k n-2 g^2 n+4 e l^2) v_{1}^2 v_{3}^3+(-2 g f o-l f j+2 g j^2-3 e m j+3 e h o+2 d j n-2 h i j-h f n-2 d i o+4 f i m-c m n+c l o) v_{1} v_{2}^4+(e i m-5 l f i-2 c l n+3 c m^2-10 e h n-b l o+k f j+6 h i^2+9 g f n+2 d i n-k c o-3 d m j+d h o+b m n-5 h f m+4 h^2 j+8 e l j+g e o-10 g i j) v_{1} v_{2}^3 v_{3}+(d h n-d l j-3 b m^2+2 b l n+3 k c n-3 c l m-18 h^2 i+6 g h j-11 g f m-8 e l i-6 g e n+d i m-5 k e j+19 e h m+k f i+10 h l f+12 g i^2-d g o+b k o) v_{1} v_{2}^2 v_{3}^2+(3 b l m+2 c l^2+12 h^3-14 e h l-10 g h i-8 d h m-4 k h f-2 g^2 j+3 d l i+4 g l f+2 d k j-2 k c m+3 d g n+7 g e m-3 b k n+7 k e i) v_{1} v_{2} v_{3}^3+(4 g^2 i+2 b k m-4 g h^2+6 d h l-4 d k i-2 b l^2+2 k e h-2 g e l-2 d g m) v_{1} v_{3}^4+(h f j-d j^2-c h o-f^2 m+d f o+c m j) v_{2}^5+(-g f j-2 d e o-4 d f n+4 c h n+2 l f^2+6 d i j+4 e f m-3 h f i-2 c l j+b h o-2 e h j-b m j-3 c i m+g c o) v_{2}^4 v_{3}+(-5 d h j-4 e^2 m+3 g f i+2 b l j+2 g e j-b g o+4 d f m+3 h^2 f+3 b i m+6 e h i-3 c h m-8 e l f+8 d e n-9 d i^2-k f^2+6 c l i-4 b h n+d^2 o+k c j-4 g c n) v_{2}^3 v_{3}^2+(-2 c h l+d g j+8 e^2 l+4 k e f-8 d e m+3 b h m-3 k c i+15 d h i-6 b l i+5 g c m-6 e h^2+4 b g n-b k j-4 d^2 n-6 g e i-4 g h f) v_{2}^2 v_{3}^3+(g^2 f+2 b h l-3 d g i-5 b g m-2 g c l+3 b k i-d k f+8 g e h-4 k e^2+2 k c h-6 d h^2+5 d^2 m) v_{2} v_{3}^4+(2 d k e-2 b k h-2 d^2 l+2 d g h-2 g^2 e+2 b g l) v_{3}^5
$
\end{center}
}}

{\fontsize{8pt}{0pt}{
\begin{center}$
\beta_{13} = (2 k i n+2 g l o-2 k h o-2 l i m+2 h m^2-2 g m n) v_{1}^4 v_{2}+(-4 k i m+4 l^2 i-4 g l n+4 g m^2-4 h l m+4 k h n) v_{1}^4 v_{3}+(3 k e o-2 k i j+5 d m n+4 l i^2-5 d l o-2 h i m+2 g h o-4 g i n-k f n-3 e m^2+2 g m j+l f m) v_{1}^3 v_{2}^2+(-2 g^2 o+4 h^2 m-4 g l j+8 g i m+10 e l m+2 k f m-2 l^2 f-12 h l i+4 k h j+10 d l n+2 d k o+2 k i^2-10 k e n-12 d m^2) v_{1}^3 v_{2} v_{3}+(8 h^2 l-8 e l^2+4 g^2 n-4 d k n-12 g h m+4 d l m-4 k h i+8 k e m+4 g l i) v_{1}^3 v_{3}^2+(4 b l o-4 l f i-4 b m n+3 d h o-5 d m j+4 g i j+h f m+c m^2+2 d i n+6 e i m+k f j-4 h i^2-6 g e o-k c o+2 g f n) v_{1}^2 v_{2}^3+(-8 g i^2-14 e h m+3 d i m+k f i+4 k c n-3 b k o-16 d h n-2 k e j+20 g e n+3 d g o-4 e l i+16 h^2 i+6 h l f+10 d l j-8 b l n-8 g h j-4 c l m+11 b m^2-7 g f m) v_{1}^2 v_{2}^2 v_{3}+(-4 k h f+4 k e i+4 g l f+4 g^2 j+12 e h l-16 h^3-4 d k j+12 g h i-8 d g n+8 b k n-16 d l i+28 d h m+4 c l^2-4 k c m-8 b l m-16 g e m) v_{1}^2 v_{2} v_{3}^2+(4 b l^2-8 g^2 i+4 d g m+8 g h^2+4 g e l+8 d k i-4 b k m-12 d h l-4 k e h) v_{1}^2 v_{3}^3+(-4 b h o-2 d i j+a m n-2 c i m+4 h f i+4 b m j+3 d e o+2 g c o+l f^2-a l o-3 e f m-d f n-2 g f j) v_{1} v_{2}^4+(6 d f m-10 b i m+6 d i^2+2 b g o+2 a l n+2 e l f-3 a m^2+6 e^2 m-3 d^2 o+a k o-k f^2-8 g c n+4 d h j+16 b h n-8 e h i+4 c l i-10 d e n+4 g e j+2 g f i-8 h^2 f-8 b l j+4 c h m) v_{1} v_{2}^3 v_{3}+(-4 d g j+11 d^2 n-8 b g n+20 b l i+10 g c m-4 g e i-14 d h i-8 c h l-2 k c i-8 e^2 l+4 b k j-16 b h m-3 a k n-7 d l f+16 e h^2+k e f+3 a l m+6 g h f+3 d e m) v_{1} v_{2}^2 v_{3}^2+(4 k c h+10 d g i-4 g c l+4 d h^2-10 b k i-12 d^2 m+2 k e^2+2 a k m+8 d e l+10 b g m+2 d k f-2 a l^2-2 g^2 f-12 g e h) v_{1} v_{2} v_{3}^3+(4 g^2 e-4 d g h-4 b g l+4 b k h-4 d k e+4 d^2 l) v_{1} v_{3}^4+(c f m-a m j-h f^2-d c o+a h o+d f j) v_{2}^5+(-2 d e j+4 e h f+4 d c n+b d o-b f m-a g o-3 d f i+g f^2-4 a h n-2 c e m-2 c l f+2 a l j+3 a i m) v_{2}^4 v_{3}+(-4 g e f-4 e^2 h+d h f-6 a l i+2 b l f+d^2 j-a k j-4 b d n+4 a g n+k c f+2 b e m+3 a h m+4 c e l+6 d e i-5 d c m) v_{2}^3 v_{3}^2+(-5 a g m+d g f+2 d c l+4 g e^2-b k f+5 b d m+3 a k i-2 k c e-4 b e l+2 a h l-2 d e h-3 d^2 i) v_{2}^2 v_{3}^3+(-2 a k h-2 b d l+2 a g l+2 b k e-2 d g e+2 d^2 h) v_{2} v_{3}^4
$
\end{center}
}}

{\fontsize{8pt}{0pt}{
\begin{center}$
\beta_{14} = (-2 l^2 i+2 g l n+2 k i m-2 g m^2+2 h l m-2 k h n) v_{1}^5+(-k f m+3 k e n-6 h^2 m-4 k i^2+5 d m^2-5 d l n+2 g h n-2 g l j-3 e l m+l^2 f+8 h l i+2 k h j) v_{1}^4 v_{2}+(-2 d l m-2 g l i+6 g h m-4 h^2 l-4 k e m+2 k h i+2 d k n-2 g^2 n+4 e l^2) v_{1}^4 v_{3}+(4 b l n+4 k f i-6 h^2 i+3 d h n-3 k e j-4 h l f-6 g e n+15 e h m-6 e l i-2 g h j+8 g i^2-4 b m^2-8 d i m-k c n+c l m+5 d l j) v_{1}^3 v_{2}^2+(-3 b k n+2 g^2 j-10 e h l+12 h^3-14 g h i+7 k e i-2 d k j+3 g e m-4 k h f-2 c l^2+7 d l i+3 b l m+4 g l f+2 k c m-8 d h m+3 d g n) v_{1}^3 v_{2} v_{3}+(4 g^2 i+2 b k m-4 g h^2+6 d h l-4 d k i-2 b l^2+2 k e h-2 g e l-2 d g m) v_{1}^3 v_{3}^2+(3 e l f+3 d e n+k c j+8 b i m+6 g e j-4 d i^2+2 c l i-a l n-3 d h j-k f^2+a m^2+2 g c n-4 b l j-8 g f i+6 e h i-9 e^2 m+4 d f m-5 c h m+3 h^2 f-4 b h n) v_{1}^2 v_{2}^3+(3 b k j-3 d g j-6 b l i+b h m-5 k c i-a l m-18 e h^2+12 e^2 l+a k n+k e f+19 d h i-11 d l f+2 b g n-g c m+d e m-3 d^2 n+6 c h l-8 g e i+10 g h f) v_{1}^2 v_{2}^2 v_{3}+(16 g e h+2 b h l+a l^2-2 g c l-2 d e l-5 g^2 f-a k m-10 d h^2+2 k c h+4 d^2 m-3 b g m+b k i-7 k e^2+5 d k f-d g i) v_{1}^2 v_{2} v_{3}^2+(2 d k e-2 b k h-2 d^2 l+2 d g h-2 g^2 e+2 b g l) v_{1}^2 v_{3}^3+(-3 e h f+2 g f^2-2 c h i+a h n+6 c e m-3 d e j-4 b f m-d c n-2 a i m-2 g c j+4 d f i-c l f+4 b h j+a l j) v_{1} v_{2}^4+(-5 g e f+8 g c i-3 d c m+3 d^2 j-10 b h i+a l i+9 b l f+a h m-a k j+d e i-2 b g j+k c f+2 b e m+6 e^2 h-a g n+4 c h^2-5 d h f+b d n-10 c e l) v_{1} v_{2}^3 v_{3}+(-12 g c h-2 a h l-2 d e h-5 d^2 i-6 b e l+4 b g i+a g m-5 b k f-b d m+8 b h^2+5 d g f+8 d c l+4 k c e+2 g e^2+a k i) v_{1} v_{2}^2 v_{3}^2+(-4 d g e+4 b k e+4 g^2 c-4 b g h-4 d k c+4 d^2 h) v_{1} v_{2} v_{3}^3+(-a h j-c^2 m+c h f+a f m-d f^2+d c j) v_{2}^5+(-b h f-2 a l f+4 d e f-3 d c i+2 b c m-b d j+a g j-g c f+2 c^2 l-2 c e h+3 a h i-2 a e m) v_{2}^4 v_{3}+(-3 a h^2+2 b e h+a k f+a d m-4 d e^2+4 d c h-2 d^2 f+2 g c e-b^2 m+b g f+4 a e l-3 a g i+3 b d i-4 b c l-k c^2) v_{2}^3 v_{3}^2+(-2 a k e+2 b k c-2 b g e-2 a d l+4 d^2 e+4 a g h-4 b d h+2 b^2 l-2 d g c) v_{2}^2 v_{3}^3+(-k b^2+k a d-g^2 a+2 b d g-d^3) v_{2} v_{3}^4
$
\end{center}
}}

{\fontsize{8pt}{0pt}{
\begin{center}$
\beta_{15} = (k f m-d m^2+e l m+d l n-k e n-l^2 f) v_{1}^5+(-4 g f m-3 d h n+6 h l f-c l m+k c n+4 g e n-3 e h m-d l j+2 b m^2-2 k f i-2 b l n+k e j-2 e l i+4 d i m) v_{1}^4 v_{2}+(-2 k h f-d g n+b k n+2 c l^2-b l m-3 d l i+2 g l f-2 e h l+3 k e i-g e m+4 d h m-2 k c m) v_{1}^4 v_{3}+(4 d f m-3 d e n+3 d h j-8 b i m+8 g f i-9 h^2 f+a l n+k f^2+2 b l j-4 d i^2+6 b h n-4 g c n-k c j+3 e^2 m+6 e h i-a m^2+2 c l i-4 g e j+3 c h m-5 e l f) v_{1}^3 v_{2}^2+(6 e h^2+a l m-b k j-3 d l f+d h i+3 d^2 n-10 g e i-a k n+4 e^2 l+9 g c m+d g j+8 b l i-2 b g n-5 b h m-5 d e m+k c i+2 g h f-10 c h l+k e f) v_{1}^3 v_{2} v_{3}+(-a l^2-3 b k i+3 d g i+4 k c h+a k m+d k f+b g m-2 d^2 m-3 k e^2+2 b h l+2 g e h-4 d h^2-g^2 f+4 d e l-4 g c l) v_{1}^3 v_{3}^2+(15 e h f-6 b h j-3 a h n-6 e^2 i+3 d e j+4 a i m+8 b i^2-2 b e n-a l j+4 g c j-4 g f^2-4 c e m-6 c h i+5 d c n-8 d f i+c l f) v_{1}^2 v_{2}^3+(-3 b d n+g e f+12 c h^2-8 b h i-6 g c i+6 c e l-b l f+19 d e i-k c f-3 d^2 j-18 e^2 h+a h m-5 a l i+10 b e m+d h f+a k j+2 b g j-11 d c m+3 a g n) v_{1}^2 v_{2}^2 v_{3}+(4 b g i-5 d^2 i-2 d e h+8 g e^2-d g f+4 a h l-6 g c h+b k f+a k i+8 d c l-12 b e l+5 b d m+2 b h^2-5 a g m-2 k c e) v_{1}^2 v_{2} v_{3}^2+(4 d^2 h-2 b g h-2 d k c+2 g^2 c-2 a k h-2 b d l-4 d g e+4 b k e+2 a g l) v_{1}^2 v_{3}^3+(-2 b c n-3 c h f+8 c e i-4 a i^2+c^2 m+2 a e n-5 d c j+5 d f^2-6 e^2 f-a f m+2 b e j+3 a h j) v_{1} v_{2}^4+(-14 b e i+3 g c f+4 b c m-4 a e m-2 a d n-3 a g j+7 d c i+3 b d j+2 b^2 n-10 c e h-8 d e f+3 b h f+2 a l f+7 a h i-2 c^2 l+12 e^3) v_{1} v_{2}^3 v_{3}+(-5 b^2 m+k c^2-7 a h^2-2 b c l-a k f-2 d c h-3 b g f+4 d^2 f-10 d e^2+16 b e h+2 a e l-b d i+2 g c e+a g i+5 a d m) v_{1} v_{2}^2 v_{3}^2+(4 b^2 l-4 b g e-4 a d l-4 b d h+4 d^2 e+4 a g h) v_{1} v_{2} v_{3}^3+(-k b^2+k a d-g^2 a+2 b d g-d^3) v_{1} v_{3}^4+(2 b c j+2 a f i-2 b f^2-2 c^2 i+2 c f e-2 a e j) v_{2}^5+(2 a d j-4 a f h-2 b^2 j+2 a e i+4 c^2 h-2 c d f-2 b c i+6 b e f-4 e^2 c) v_{2}^4 v_{3}+(-2 b c h-2 b d f+6 d c e+2 a g f-4 b e^2-2 g c^2+2 a e h+4 b^2 i-4 a d i) v_{2}^3 v_{3}^2+(2 b g c-2 d^2 c+2 b d e-2 a g e-2 b^2 h+2 a d h) v_{2}^2 v_{3}^3
$
\end{center}
}}

{\fontsize{8pt}{0pt}{
\begin{center}$
\beta_{16} = (-d m n+d l o-l f m+e m^2-k e o+k f n) v_{1}^5+(2 l f i-c m^2+2 b m n+3 h f m+d m j+2 d i n-4 e i m-3 d h o-k f j-4 g f n+4 g e o-2 b l o+k c o) v_{1}^4 v_{2}+(-4 d l j+2 h l f-4 e h m-3 k f i-b m^2+3 d i m-2 k c n+g f m+2 d h n-d g o+2 c l m+4 k e j+b k o) v_{1}^4 v_{3}+(5 d f n-4 g c o+a l o-2 d i j-a m n-l f^2+4 c i m-6 h f i-e f m+4 e i^2-4 b i n+6 b h o+4 g f j-2 b m j-3 d e o) v_{1}^3 v_{2}^2+(-4 e l f-6 h^2 f-2 c h m-4 c l i+10 d h j-2 b g o+8 g c n-a k o+10 g f i+3 d^2 o-6 d i^2-2 k c j+3 k f^2+8 e^2 m-16 g e j-4 b h n-2 b i m+8 b l j-6 d f m+a m^2+8 e h i-4 d e n) v_{1}^3 v_{2} v_{3}+(-4 b k j-a l m+6 k c i-2 g c m-2 g h f-4 c h l+a k n-d e m+5 d l f+4 e h^2+4 d g j+4 b h m-6 d h i-3 k e f-d^2 n) v_{1}^3 v_{3}^2+(-2 b e o+5 d c o-c f m+4 b i j-4 c i^2+3 h f^2+2 a i n+a m j-2 b f n-3 a h o-5 d f j+2 e f i) v_{1}^2 v_{2}^3+(-16 e^2 i+2 a h n-a i m+3 a g o+16 d e j-6 c e m+8 b i^2-3 d f i-3 b d o+8 g c j-10 d c n+4 c l f-11 g f^2+8 b e n-4 a l j-20 b h j+7 b f m+14 e h f+4 c h i) v_{1}^2 v_{2}^2 v_{3}+(-a h m-20 g c i+8 b g j+2 a l i+3 a k j-3 k c f+4 b d n+8 c e l-4 a g n+16 g e f+14 d e i-6 b e m-10 b l f+7 d c m-16 e^2 h-11 d^2 j-3 d h f+4 b h i+8 c h^2) v_{1}^2 v_{2} v_{3}^2+(5 b k f-3 a k i-5 d g f-4 b h^2-2 k c e-2 d c l+2 d e h+3 d^2 i+a g m+4 g c h-b d m+2 a h l) v_{1}^2 v_{3}^3+(-2 e f^2-2 a i j+2 b f j-2 b c o+2 a e o+2 c f i) v_{1} v_{2}^4+(-4 a e n-10 d c j+2 c^2 m-4 e^2 f+10 a h j-2 a d o+12 c e i-2 a f m+4 b c n+2 b^2 o-2 a i^2-8 b f i+12 d f^2-10 c h f) v_{1} v_{2}^3 v_{3}+(16 e^3+8 b d j+4 a d n+4 a e m-4 c^2 l+4 a l f+16 b h f-12 c e h-8 a g j-4 b^2 n+8 g c f+16 d c i-4 a h i-28 d e f-4 b c m-12 b e i) v_{1} v_{2}^2 v_{3}^2+(12 d^2 f-4 d e^2-2 a k f+10 a g i-2 a d m+12 b e h-4 a e l-10 b d i-10 b g f+2 b^2 m-2 a h^2+2 k c^2-8 d c h+4 b c l) v_{1} v_{2} v_{3}^3+(2 d g c-2 d^2 e-2 b k c+2 b d h-2 a g h+2 a k e) v_{1} v_{3}^4+(-4 b f^2+4 a f i+4 c f e-4 a e j+4 b c j-4 c^2 i) v_{2}^4 v_{3}+(-8 a f h-4 c d f+4 a e i+12 b e f-4 b^2 j+8 c^2 h-8 e^2 c-4 b c i+4 a d j) v_{2}^3 v_{3}^2+(-8 a d i+12 d c e+4 a g f-8 b e^2-4 g c^2+4 a e h-4 b c h-4 b d f+8 b^2 i) v_{2}^2 v_{3}^3+(4 b d e+4 a d h-4 b^2 h+4 b g c-4 a g e-4 d^2 c) v_{2} v_{3}^4
$
\end{center}
}}

{\fontsize{8pt}{0pt}{
\begin{center}$
\beta_{17} = (-k h o+g l o-l i m-g m n+k i n+h m^2) v_{1}^5+(-2 g i n+k e o-h i m+g m j+g h o-k i j+2 d m n-2 d l o+2 l i^2-e m^2) v_{1}^4 v_{2}+(-4 g l j+l^2 f+e l m-g^2 o+4 k h j-d l n-4 h^2 m+d k o+4 g i m+2 h l i-3 k i^2-k e n-k f m+2 g h n) v_{1}^4 v_{3}+(-b m n+2 g i j+b l o+2 d i n-g f n-2 h i^2+e i m-2 d m j+2 h f m-g e o-l f i) v_{1}^3 v_{2}^2+(-9 d i m-3 k e j-6 h l f-6 g h j-d h n+9 e h m+2 h^2 i+g f m+9 d l j-b k o+5 k f i+2 b l n+d g o-6 e l i+4 g i^2-b m^2+c l m-k c n) v_{1}^3 v_{2} v_{3}+(4 g l f+2 k c m-2 c l^2+2 d l i+4 g^2 j+2 d h m-4 k h f-4 g e m-8 g h i-2 e h l-4 d k j+6 k e i+4 h^3) v_{1}^3 v_{3}^2+(c i m-2 d i j-h f i-g c o+2 e i^2+g f j-2 e f m+b h o+2 d f n-2 b i n+b m j) v_{1}^2 v_{2}^3+(4 g e j+a m^2-10 g f i-a l n-6 b l j+6 g c n-d e n+10 b i m+7 e l f-k f^2+k c j+5 h^2 f-5 e^2 m+d h j-4 b h n-7 c h m) v_{1}^2 v_{2}^2 v_{3}+(-4 k c i+a k n+4 e^2 l-k e f+7 d e m-d^2 n-14 e h^2+10 c h l-4 b l i+4 b k j-4 d g j-9 d l f+8 d h i-a l m+10 g h f-6 g c m) v_{1}^2 v_{2} v_{3}^2+(8 g e h-3 k e^2+5 d k f-2 d e l+a l^2-5 g^2 f-a k m+b g m-4 d h^2-b k i+d g i) v_{1}^2 v_{3}^3+(2 d c o+2 b i j-2 c i^2-b e o-2 d f j-b f n-a h o+a i n+h f^2+e f i) v_{1} v_{2}^4+(-2 g c j+6 c h i-4 b i^2+6 c e m+d e j-b f m-9 d c n+9 d f i+6 b e n-9 e h f-c l f-5 a i m+a g o-b d o+g f^2+a l j-2 e^2 i+3 a h n) v_{1} v_{2}^3 v_{3}+(-7 d h f-4 a g n-10 c e l-8 d e i+k c f-10 b e m+a h m+4 b d n+d^2 j+14 e^2 h+9 d c m-a k j+4 g c i+6 b l f+4 a l i-4 c h^2) v_{1} v_{2}^2 v_{3}^2+(4 b e l-4 a h l-4 b k f+4 b h^2+4 a g m-4 g e^2-4 g c h+4 k c e+4 d g f-4 b d m) v_{1} v_{2} v_{3}^3+(b d l+b k e+a k h+d^2 h-2 b g h-d g e-2 d k c+2 g^2 c-a g l) v_{1} v_{3}^4+(-a i j+c f i+a e o+b f j-b c o-e f^2) v_{2}^5+(d c j-c h f+b^2 o+3 a i^2-4 a e n+4 e^2 f-c^2 m-2 b e j-2 c e i+a h j-a d o-4 b f i+4 b c n+a f m) v_{2}^4 v_{3}+(-4 e^3+4 a e m-2 d e f+2 c e h-2 a l f-2 d c i+2 c^2 l+4 a d n-6 a h i+8 b e i+4 b h f-4 b^2 n-4 b c m) v_{2}^3 v_{3}^2+(2 d c h-5 a d m-b d i+a g i+5 b^2 m-8 b e h+a k f+4 d e^2+3 a h^2-k c^2-b g f) v_{2}^2 v_{3}^3+(2 b g e-a k e-d^2 e+b k c+2 a d l-a g h-d g c+b d h-2 b^2 l) v_{2} v_{3}^4
$
\end{center}
}}

{\fontsize{8pt}{0pt}{
\begin{center}$
\beta_{18} = (2 h l n-g m n-h m^2+l i m-2 l^2 j-k i n+2 k m j-k h o+g l o) v_{1}^5+(-4 e l n-2 l i^2+3 d m n+k f n+h i m-l f m-7 g m j+10 h l j+6 g i n-3 d l o+g h o-6 h^2 n-3 k i j+2 k e o+2 e m^2) v_{1}^4 v_{2}+(-d m^2+4 l^2 f+3 k i^2+4 h^2 m+d k o-g^2 o+2 g i m-4 k f m-8 h l i) v_{1}^4 v_{3}+(-c m^2-5 g e o-12 h^2 j-k c o-12 d i n+3 l f i+h f m-8 e l j+9 d m j-3 b m n-5 g f n-5 e i m+2 c l n+10 g i j+18 e h n+3 b l o+2 h i^2+k f j+3 d h o) v_{1}^3 v_{2}^2+(-16 h l f+16 h^2 i-2 d i m-2 b k o+2 b m^2-16 g i^2+14 g f m-16 e h m+16 e l i+2 k f i+2 d g o) v_{1}^3 v_{2} v_{3}+(6 e h l-2 d k j-3 k e i+2 g^2 j-4 h^3-2 g l f+b k n-b l m+d l i-3 g e m+2 d h m-2 c l^2-d g n+2 k h f+2 k c m+2 g h i) v_{1}^3 v_{3}^2+(-5 h f i+9 d f n-5 b m j-8 c h n+2 c l j+18 e h j+2 e i^2-3 g f j+10 b i n+e f m-12 d i j-5 b h o+3 g c o+a m n-a l o+3 c i m-12 e^2 n+3 d e o-l f^2) v_{1}^2 v_{2}^3+(2 b g o-3 d^2 o-a m^2-2 b i m-40 e h i+a k o+8 e l f-20 d f m+16 e^2 m+16 h^2 f-8 c l i+26 d i^2+8 c h m-2 g f i-k f^2) v_{1}^2 v_{2}^2 v_{3}+(-k e f+5 d l f-3 d g j-k c i-12 e^2 l-b h m-a k n+3 d^2 n-2 b g n-5 g c m+6 c h l+14 g e i+6 e h^2+a l m-13 d h i+3 b k j+5 d e m-4 g h f) v_{1}^2 v_{2} v_{3}^2+(-2 d^2 m-2 g^2 f-2 b k i+2 b g m+2 d g i+2 d k f) v_{1}^2 v_{3}^3+(-c f m+3 d f j+2 h f^2-3 d c o+6 b i j+a m j+b e o-4 c h j-2 c i^2-6 e^2 j-3 a i n+10 c e n+2 a h o-7 b f n+e f i) v_{1} v_{2}^4+(-2 d f i+16 c h i-16 e h f+16 e^2 i+2 a i m-16 c e m+2 g f^2+14 b f m-16 b i^2-2 a g o+2 b d o) v_{1} v_{2}^3 v_{3}+(5 d h f-a k j-12 c h^2-4 b e m+k c f-a h m+6 c e l-3 b d n+14 b h i-13 d e i-2 b g j-a l i+3 a g n-g e f+5 d c m+6 e^2 h+3 d^2 j-5 b l f) v_{1} v_{2}^2 v_{3}^2+(2 b d m+2 d g f-2 d^2 i-2 a g m-2 b k f+2 a k i) v_{1} v_{2} v_{3}^3+(-a k h+2 g^2 c-2 b g h+3 b k e+a g l+3 d^2 h-2 d k c-b d l-3 d g e) v_{1} v_{3}^4+(b c o-2 c^2 n+2 a f n-a e o-e f^2-a i j+c f i+2 c e j-b f j) v_{2}^5+(-b^2 o+4 c^2 m+a d o-d f^2+4 e^2 f+2 b f i+3 a i^2-4 a f m-8 c e i) v_{2}^4 v_{3}+(-2 b c m-2 c^2 l+d c i+2 a l f-4 e^3-3 b h f-g c f+a g j+2 b e i+2 d e f-2 a d n+6 c e h+2 a e m+2 b^2 n-3 a h i-b d j) v_{2}^3 v_{3}^2+(2 a d m-2 d^2 f+2 b d i-2 b^2 m+2 b g f-2 a g i) v_{2}^2 v_{3}^3+(-3 b d h-2 a d l-d g c+3 a g h+3 d^2 e+2 b^2 l+b k c-2 b g e-a k e) v_{2} v_{3}^4+(-k b^2+k a d-g^2 a+2 b d g-d^3) v_{3}^5
$
\end{center}
}}
The coefficients of these polynomials are the explicit matrix elements of the $21 \times 21$ undulation matrix. The MAPLE program that was used to generate these explicit expressions is available upon request.


\begin{thebibliography}{100}

\bibitem{CayleySalmon} A. Cayley and G. Salmon, \emph{A Treatise on the Higher Plane Curves}, Elibron Classics, original 1852

\bibitem{Mandelbrot}
    V.Dolotin and A.Morozov, \emph{Algebraic Geometry of Discrete Dynamics. The case of one variable}, hep-th/0501235; \emph{The Universal Mandelbrot Set: Beginning of the Story}, World Scientific, 2006 \\
    J.Milnor, \emph{Dynamics of one complex variable}, Princeton University Press, 2006 \\
    An.Morozov, \emph{Universal Mandelbrot Set as a Model of Phase Transition Theory}, JETP Lett. \textbf{86} (2007) 745-748, arXiv:0710.2315

\bibitem{RiemannSchottky} I.Krichever and T. Shiota, \emph{Soliton equations and the Riemann-Schottky problem}, arXiv:1111.0164

\bibitem{Nolinal} V.Dolotin and A.Morozov, \emph{Introduction to Non-linear Algebra}, World Scientific, 2007; hep-th/0609022

\end{thebibliography}
\end{document}